\documentclass[12pt]{article}

\usepackage{amssymb}
\usepackage{mathrsfs}
\usepackage{theorem}
\def\math{\mathscr}

\newtheorem{theorem}{Theorem}
\newtheorem{lemma}{Lemma}
\newtheorem{corollary}{Corollary}
\newenvironment{proof} {{\bf Proof.}}{\hfill \fbox{}\\ \smallskip}

\newcommand{\paragrafo}[1]{\section{#1}\setcounter{equation}{0}}

{\theorembodyfont{\rmfamily} \newtheorem{remark}{Remark}}
{\theorembodyfont{\rmfamily} \newtheorem{definition}{Definition}}
{\theorembodyfont{\rmfamily} \newtheorem{example}{Example}}

\newcommand{\reff}[1]{{\rm(\ref{#1})}}

\def\display{\displaylines}
\def\grande{\displaystyle}

\newcommand{\beginrighe}[1]{\begin{eqnarray}\label{#1}\begin{array}{c}\grande}
\def\endrighe{\end{array}\end{eqnarray}}
\def\crg{\cr\grande}

\def\C{\mathbb{C}}
\def\R{\mathbb{R}}
\def\N{\mathbb{N}}

\def\de{\partial}

\def\a{\alpha}
\def\b{\beta}
\def\g{\gamma}

\def\D{\Delta}
\def\e{\varepsilon}

\def\h{\eta}
\def\l{\lambda}

\def\m{\mu}

\def\r{\varrho}
\def\s{\sigma}

\def\f{\varphi}

\def\O{\Omega}
\def\o{\omega}
\def\p{\psi}

\def\Cspt{C_{0}}
\def\Hspt{H_{0}}
\def\Wspt{W_{0}}

\def\Re{\mathop{\mathscr Re}\nolimits}
\def\Im{\mathop{\mathscr Im}\nolimits}
\def\elle{\mathop{\mathscr L}\nolimits}
\def\A{\mathop{\mathscr A}\nolimits}
\def\B{\mathop{\mathscr B}\nolimits}
\def\Sm{\mathop{\mathscr S}\nolimits}
\def\Cm{\mathop{\mathscr C}\nolimits}

\def\Dom{{\math D}}
\def\spt{\mathop{\rm spt}\nolimits}

\def\dive{\nabla^{t}}

\def\leq{\leqslant}
\def\geq{\geqslant}
\def\lan{\langle}
\def\ran{\rangle}

\title{
    Criterion for the $L^{p}$-dissipativity of second order differential
    operators with complex coefficients
}
\author{A. Cialdea 
\thanks{Dipartimento di Matematica, Universit\`a della Basilicata, Viale 
dell'Ateneo Lucano 10, 85100, Potenza, Italy. \textit{email:}
cialdea@email.it.}  
\and V. Maz'ya
\thanks{Department of Mathematics, Ohio State University,
231 W 18th Avenue, Columbus, OH 43210, USA.
Department of Mathematical Sciences, M{\&}O Building,
University of Liverpool, Liverpool L69 3BX, UK.
\textit{email:} vlmaz@mai.liu.se.}
}
\date{}    

\begin{document}

    \maketitle

   {\small  {\bf Abstract.}
    We prove that the algebraic condition 
    $|p-2|\, |\lan \Im\A\xi,\xi\ran| \leq 2 \sqrt{p-1}\,
	    \lan \Re\A\xi,\xi\ran$ (for any $\xi\in\R^{n}$) is
    necessary and
    sufficient for the $L^{p}$-dis\-sipativity
    of the Dirichlet problem for
    the differential
    operator $\nabla^{t}(\A\nabla)$, where $\A$ is a matrix whose
    entries are complex measures and whose imaginary part is symmetric.
    This result is new even for
    smooth coefficients, when it implies a criterion for the
    $L^{p}$-contractivity of the corresponding semigroup.
    We consider also the operator
    $\nabla^{t}(\A\nabla)+{\bf b}\nabla +a$, where the coefficients
    are smooth and $\Im\A$ may be not symmetric.
    We show that the previous algebraic
    condition is necessary and sufficient for the
    $L^{p}$-quasi-dissipativity of this operator. The same
    condition is necessary and sufficient 
    for the $L^{p}$-quasi-contractivity of the corresponding 
    semigroup.
    We give a necessary and sufficient condition for the 
    $L^{p}$-dissipativity in $\R^{n}$ of the operator
    $\nabla^{t}(\A\nabla)+{\bf b}\nabla +a$ with constant coefficients.}
    \smallskip
    
    {\small {\bf R\'esum\'e.}
    On montre que la condition alg\'ebrique
    $|p-2|\, |\lan \Im\A\xi,\xi\ran| \leq 2 \sqrt{p-1}\,
		\lan \Re\A\xi,\xi\ran$ (pour tout $\xi\in\R^n$) est
    n\'ecessaire et suffisante pour la dissipativit\'e $L^p$
    du probl\`eme de Dirichlet pour l'op\'erateur diff\'erentiel
    $\nabla^{t}(\A\nabla)$, o\`u $\A$ est une matrice
    dont les coefficients sont des mesures complexes et dont la partie
    imaginaire est sym\'etrique. Ce r\'esultat est nouveau m\^eme pour
    des coefficients r\'eguliers, quand il implique un crit\`ere pour
    la contractivit\'e $L^p$ du semi-groupe correspondant. On consid\`ere
    aussi l'op\'erateur $\nabla^{t}(\A\nabla)+{\bf b}\nabla +a$,
    o\`u les coefficients sont r\'eguliers et $\Im\A$ n'est pas
    n\'ecessairement sym\'etrique. On montre que la condition alg\'ebrique
    pr\'ec\'edente est n\'ecessaire et suffisante pour la
    quasi-dissipativit\'e $L^p$ de cet op\'erateur. La m\^eme condition est
    n\'ecessaire et suffisante pour la quasi-contractivit\'e $L^p$ du
    semi-groupe correspondant. On donne une condition n\'ecessaire et
    suffisante pour la dissipativit\'e $L^p$ dans $\R^n$ de l'op\'erateur
    $\nabla^{t}(\A\nabla)+{\bf b}\nabla +a$ avec des coefficients
    constants.}
    
    \medskip
    
\paragrafo{Introduction}

Various aspects of the
$L^{p}$-theory of semigroups generated by linear differential operators
were studied in
\cite{brezis, davies, amann, strichartz2,
davies1, kovalenko, robinson, davies2, davies3,
liskevich1, ouha, langermazya, langer,
daners, karrmann, sobol, liskevich2} \textit{et al}.
In particular, it has been known for years that scalar second order
elliptic operators with
real coefficients may generate contractive semigroups in $L^{p}$
\cite{mazyasobolev}.

Necessary and sufficient conditions for the
$L^{\infty}$-contractivity for general second order strongly elliptic systems
with smooth coefficients
were given in \cite{kresin}, where
 scalar second order elliptic operators \textit{with
complex coefficients} were handled as  a particular case.
Such operators  generating
$L^{\infty}$-contractive semigroups were later characterized in \cite{auscher}
under the assumption that the coefficients are measurable and bounded.

In the present paper we find an algebraic necessary and sufficient
condition for the $L^{p}$-dissipativity of the Dirichlet problem for
the differential operator
$$
    A= \nabla^{t}(\A\nabla)
$$
where $\A$ is a matrix whose
entries are complex measures and whose imaginary part is symmetric. 
Namely in Section \ref{sec:main}, after giving the definition of 
$L^{p}$-dissipativity of the corresponding form  
$$
{\elle}(u,v)=\int_{\O}\lan \A\nabla u,\nabla v\ran \ ,
$$
we prove that
$\elle$ is $L^{p}$-dissipative if and only if
\begin{equation}
    |p-2|\, |\lan \Im\A\xi,\xi\ran| \leq 2 \sqrt{p-1}\,
        \lan \Re\A\xi,\xi\ran
    \label{eq:0}
\end{equation}
for any $\xi\in\R^{n}$.   This result is new even for
smooth coefficients. An example shows that the
statement is not true if $\Im\A$
is not symmetric.

It is impossible, in general,  to obtain a similar algebraic
characterization  for the operator with lower order terms
\begin{equation}
    Au =\dive(\A\nabla u) + {\bf b}\nabla u + \dive({\bf c}u)+au.
    \label{eq:N0}
\end{equation}

In fact, consider
for example the operator
$$
Au=\D u + a(x)u
$$
in a bounded domain $\O\subset\R^{n}$.  Denote by $\l_{1}$ 
the first eigenvalue of the Dirichlet problem for Laplace equation 
in $\O$. A sufficient condition for $A$ 
to be $L^{2}$-dissipative is $\Re a\leq \l_{1}$ and we cannot
give an algebraic characterization of 
$\l_{1}$.
However in Section \ref{sec:constant} we give a 
necessary and sufficient condition
for the $L^{p}$-dissipativity of 
operator \reff{eq:N0} in $\R^{n}$ for
the particular case of constant coefficients.

In  Section \ref{sec:smooth} 
we consider operator \reff{eq:N0} with smooth coefficients
without the requirement of simmetricity of $\Im\A$. 
After showing that
the concept of $L^{p}$-dissipativity of the form  $\elle$
 is equivalent to the usual
 $L^{p}$-dissipativity of the
operator $A$, we
prove that the algebraic condition \reff{eq:0} is, in general, 
necessary and sufficient for the $L^{p}$-quasi-dissipativity, i.e. 
for the
$L^{p}$-dissipativity of $A-\o I$ for a suitable $\o>0$.

In other words the range of the exponent $p$ admissible for the 
$L^{p}$-quasi-dissipa\-tivity is given by the inequalities
$$
              2+2\l(\l-\sqrt{\l^{2}+1}) \leq p \leq 
                   2+2\l(\l+\sqrt{\l^{2}+1}),
$$
where 
$$
\l=\inf_{(\xi,x)\in {\cal M}}{\lan \Re\A(x)\xi,\xi\ran \over
  |\lan \Im\A(x)\xi,\xi\ran|}
$$
and ${\cal M}=\{(\xi,x)\in\R^{n}\times \O\ |\ 
\lan \Im\A(x)\xi,\xi\ran\neq 0\}$.

Finally we show that \reff{eq:0} is necessary and sufficient for the 
$L^{p}$-quasi-contractivity of the semigroup generated by the 
Dirichlet problem for the operator \reff{eq:N0}.

\paragrafo{Preliminaries}

Let $\O$ be an open set in $\R^{n}$.
By $\Cspt(\O)$ we denote the space of  complex valued continuous functions
having compact support in $\O$. Let $\Cspt^{1}(\O)$
consist of all the functions in $\Cspt(\O)$ having continuos partial
derivatives of the first order.
The inner product either in
$\C^{n}$ or in $\C$ is denoted by $\lan \cdot, \cdot \ran$  and,
as usual, the bar denotes complex conjugation.

In what follows, $\A$ is a $n\times n$ matrix
function with complex valued entries $a^{hk}\in (\Cspt(\O))^{*}$,
$\A^{t}$ is its transposed matrix and $\A^{*}$ is its adjoint
matrix, i.e. $\A^{*}=\overline{\A}^{t}$.

Let
${\bf b}=(b_{1},\ldots,b_{n})$ and ${\bf c}=(c_{1},\ldots,c_{n})$
stand for complex
valued vectors
with $b_{j}, c_{j}\in (\Cspt(\O))^{*}$.
By $a$ we mean a complex valued scalar distribution in $(\Cspt^{1}(\O))^{*}$.

We denote by ${\mathscr L}(u,v)$ the sesquilinear form
$$
    {\elle}(u,v)=\int_{\O}(\lan \A\nabla u,\nabla v\ran -\lan {\bf b}\nabla u,v\ran
    +\lan u,\overline{\bf c}\nabla v\ran
    -a\lan u,v\ran)\,
$$
defined on $\Cspt^{1}(\O)\times \Cspt^{1}(\O)$.

If $p\in(1,\infty)$, $p'$ denotes its conjugate exponent $p/(p-1)$.

\begin{definition} Let $1<p<\infty$. The form ${\mathscr L}$ is
called \textit{$L^{p}$-dissipative} if
for all $u\in \Cspt^{1}(\O)$
\begin{eqnarray}
    \Re {\mathscr L}(u, |u|^{p-2}u) \geq 0
    \qquad \hbox{\rm if}\ p\geq 2;
    \label{eq:defdis1}\\
    \Re {\mathscr L}(|u|^{p'-2}u, u) \geq 0 \qquad
    \hbox{\rm if}\ 1<p< 2
      \label{eq:defdis2}
\end{eqnarray}
(we use here that $|u|^{q-2}u\in \Cspt^{1}(\O)$ for $q\geq 2$ and
$u\in \Cspt^{1}(\O)$).
\end{definition}

The form ${\mathscr L}$ is related to the operator
\begin{equation}
    Au =\dive(\A\nabla u) + {\bf b}\nabla u + \dive({\bf c}u)+au.
    \label{eq:N}
\end{equation}
where $\dive$ denotes the divergence operator.
The operator $A$ acts from $\Cspt^{1}(\O)$ to $(\Cspt^{1}(\O))^{*}$ through the
relation
$$
\elle(u,v)=\int_{\O}\lan Au,v\ran
$$
for any $u,v \in \Cspt^{1}(\O)$.

We start with the following Lemma

\begin{lemma}\label{lemma:1}

 The form ${\mathscr L}$ is $L^{p}$-dissipative
    if and only if for all $v\in \Cspt^{1}(\O)$
    \beginrighe{eq:22}
    \grande
        \Re \int_{\O}\Big[ \lan \A\nabla v,\nabla v\ran -
	(1-2/p)\lan (\A-\A^{*})\nabla(|v|),|v|^{-1}\overline{v}\nabla v\ran  -\\
	\grande
	(1-2/p)^{2}\lan \A \nabla(|v|),\nabla(|v|)\ran
\Big]
+
\grande \int_{\O}\lan\Im ({\bf b}+{\bf c}),\Im(\overline{v}\nabla
v)\ran\,   +
\cr
\grande
\int_{\O}\Re(\dive ({\bf b}/p - {\bf c}/p') - a
)|v|^{2}
	\geq 0.
    \endrighe
  Here and in the sequel the integrand is extended by zero
    on the set where $v$ vanishes.
\end{lemma}

\begin{proof}

    {\it Sufficiency.}
    Let us prove the sufficiency for $p\geq 2$.
Suppose \reff{eq:22} holds, take $u\in \Cspt^{1}(\O)$ and
set
$$
v=|u|^{p-2\over 2}u.
$$

Since $p\geq 2$ we have $v\in \Cspt^{1}(\O)$. Moreover,
$u=|v|^{2-p\over p}v$ and therefore
$$\display{
\lan \A\nabla u, \nabla(|u|^{p-2}u)\ran =
\lan \A\nabla(|v|^{2-p\over p}v), \nabla(|v|^{p-2\over
p}v)\ran =\cr
\left\lan \A\left(\nabla v-\right(1-2/p)|v|^{-1}v\nabla |v|),
\nabla v +(1-2/p) |v|^{-1}v \nabla |v|
\right\ran=\cr
\lan \A\nabla v,\nabla v\ran -(1-2/p)\left(
\lan |v|^{-1}v\A\nabla |v|,\nabla v\ran -
\lan \A\nabla v, |v|^{-1}v\nabla |v|\ran \right)
- \cr
-\left(1-2/p\right)^{2} \lan \A \nabla |v|, \nabla |v|\ran
}
$$
Since
$$\display{
\Re (
\lan v\A\nabla |v|,\nabla v\ran -
\lan \A\nabla v, v\nabla |v|\ran ) =\cr
 \Re (v\lan \A\nabla |v|,\nabla v\ran -
\overline{\lan v\A^{*}\nabla |v|, \nabla v\ran})=
\Re (\lan v(\A-\A^{*})\nabla |v|, \nabla v\ran)
}
$$
we have
$$\display{
\Re \lan \A\nabla u,\nabla(|u|^{p-2}u)\ran=
\Re \Big[
\lan \A\nabla v,\nabla v\ran -\cr
	(1-2/p)\lan (\A-\A^{*})\nabla(|v|),|v|^{-1}\overline{v}\nabla v\ran 
	 -
	(1-2/p)^{2}\lan \A \nabla(|v|),\nabla(|v|)\ran
\Big].
}
$$

Moreover, we have
$$
\lan {\bf b} \nabla u, |u|^{p-2}u\ran = (1-2/p)\, |v|\,
{\bf b}\nabla |v| + \overline{v}\, {\bf b}\nabla v
$$
and then
$$\display{
\Re \lan {\bf b}\nabla u, |u|^{p-2}u\ran = 2\, \Re (
{\bf b}/p)\Re(\overline{v}\nabla v) -(\Im
{\bf b})\Im(\overline{v}\nabla v)=\cr
\Re ({\bf b}/p)\nabla(|v|^{2})-(\Im
{\bf b})\Im(\overline{v}\nabla v).
}
$$

An integration by parts gives
\begin{equation}
    \int_{\O}\Re \lan {\bf b}\nabla u, |u|^{p-2}u\ran=
    -\int_{\O}\Re(\nabla^{t}({\bf b}/p))|v|^{2}-
    \int_{\O}\lan \Im{\bf b},\Im(\overline{v}\nabla v)\ran
    \, .
    \label{eq:part1}
\end{equation}

In the same way we find
$$\display{
\Re \lan u,\overline{\bf c}\nabla(|u|^{p-2}u)\ran =
\Re\left( (1-2/p)\, |v| {\bf c}\nabla|v| + v\,
{\bf c}\nabla\overline{v}\right) =\cr
2\, \Re ({\bf c}/p')\Re(\overline{v}\nabla v) +(\Im
{\bf c})\Im(\overline{v}\nabla v) =\cr
\Re ({\bf c}/p')\nabla(|v|^{2}) + (\Im
{\bf c})\Im(\overline{v}\nabla v)
}
$$
and then
\begin{equation}
    \int_{\O}\Re \lan u,\overline{\bf c}\nabla(|u|^{p-2}u)\ran=
    -\int_{\O}\Re(\nabla^{t}({\bf c}/p')|v|^{2}+
    \int_{\O}\lan \Im {\bf c}, \Im(\overline{v}\nabla v)
    \ran .
    \label{eq:part2}
\end{equation}

Finally, since we have also
$$
\Re(a\lan u,|u|^{p-2}u\ran=(\Re a)|u|^{p}=(\Re a)|v|^{2},
$$
the left-hand side in \reff{eq:22} is equal to $\Re
\elle(u,|u|^{p-2}u)$ and
\reff{eq:defdis1} follows from \reff{eq:22}.

Let us suppose that $1<p<2$. Now
\reff{eq:defdis2} can be
written as
\beginrighe{eq:M}
\grande \Re \int_{\O}( \lan \A^{*}\nabla u, \nabla(|u|^{p'-2}u)\ran +
\lan \overline{\bf c}\nabla u, |u|^{p'-2}u\ran - \lan \nabla u, {\bf
b}\nabla(|u|^{p'-2}u)\ran -\cr
\grande -a\lan u, |u|^{p'-2}u\ran)\,   \geq 0.
\endrighe

We know that this is true if
\beginrighe{eq:P}
\grande	\Re \int_{\O}\Big[ \lan \A^{*}\nabla v,\nabla v\ran -
	(1-2/p')\lan (\A^{*}-\A)\nabla(|v|),|v|^{-1}\overline{v}\nabla 
	v\ran  -
	\cr
-	(1-2/p')^{2}\lan \A^{*} \nabla(|v|),\nabla(|v|)\ran
\Big]
+\cr
\grande
+\int_{\O}\lan\Im (-\overline{\bf c}-
\overline{\bf b}),\Im(\overline{v}\nabla
v)\ran\,   +
\cr
\grande
\int_{\O}\Re\left[ \dive\left((-\overline{\bf c})/p'
- (-\overline{\bf b})/p\right)-a\right]|v|^{2}
	\geq 0
\endrighe
for any $v\in \Cspt^{1}(\O)$. This condition is exactly \reff{eq:22}
and the sufficiency is proved also for $1<p<2$.

{\it Necessity.}
Let us suppose \reff{eq:defdis1} holds. Let $v\in \Cspt^{1}(\O)$
 and set
\begin{equation}
     g_{\e}=(|v|^{2}+\e^{2})^{1\over 2},\quad u_{\e}=g_{\e}^{{2\over
    p}-1}v.
    \label{eq:defge}
\end{equation}

We have
$$\displaylines{
 \lan \A \nabla u_{\e},\nabla(|u_{\e}|^{p-2}u_{\e}) \rangle =\cr
|u_{\e}|^{p-2}\langle \A\nabla u_{\e}, \nabla u_{\e}\rangle +
(p-2)|u_{\e}|^{p-3}\langle
\A\nabla u_{\e},u_{\e}\nabla |u_{\e}|\rangle
}
$$

A direct computation shows that
       $$\displaylines{
       |u_{\e}|^{p-2} \langle \nabla u_{\e}, \nabla u_{\e}\rangle =
       \left(1-2/p\right)^{2}g_{\e}^{-(p+2)}|v|^{p+2}\lan \nabla |v|,
       \nabla |v|\ran  -
       \cr
       \left(1-2/p\right)g_{\e}^{-p}|v|^{p-1}(\lan v\nabla |v|,
       \nabla v\ran + \lan\nabla v, v\, \nabla |v|) +
       g_{\e}^{2-p}|v|^{p-2}\lan \nabla v, \nabla v\ran\, ,
       }
       $$
       $$\displaylines{
       |u_{\e}|^{p-3}\langle
       \nabla u_{\e},u_{\e}\nabla |u_{\e}|\rangle = \cr
       \left[\left(1-2/p\right)^{2}g_{\e}^{-(p+2)}|v|^{p+2}-
       \left(1-2/p\right)g_{\e}^{-p}|v|^{p}\right]
       \lan \nabla |v|, \nabla |v|\ran + \cr
       \left[-\left(1-2/p\right)g_{\e}^{-p}|v|^{p-1}+
       g_{\e}^{-p+2}|v|^{p-3}\right]
    \lan \nabla v, v \nabla |v|\ran .
       }
       $$

Observing that $g_{\e}$ tends to $|v|$ as $\e \to 0$ and referring to
Lebesgue's
dominated convergence theorem we find

\beginrighe{zac1}
\lim_{\e\to 0}\int_{\O}\langle
\A\nabla u_{\e},\nabla(|u_{\e}|^{p-2}u_{\e}) \rangle   =\crg
\int_{\O}\langle \A\nabla v, \nabla v\rangle\,   -\crg
(1-2/p)\int_{\O}{1\over |v|} \left(
\langle v\A  \nabla |v|, \nabla v\rangle-
\langle \A\nabla v,v\, \nabla |v|\rangle\right)\,   -\crg
- \left(1-2/p\right)^{2} \int_{\O}
\langle
\A\nabla |v|, \nabla |v|\rangle\,   \ .
\endrighe

Similar computations show that
$$\display{
\lan {\bf b}\nabla u_{\e},|u_{\e}|^{p-2}u_{\e}\ran =
-(1-2/p)g_{\e}^{-p}|v|^{p+1}{\bf b}\nabla |v| +
g_{\e}^{2-p}|v|^{p-2}\overline{v}{\bf b}\nabla v\cr
\lan u_{\e},\overline{\bf c}\nabla(|u_{\e}|^{p-2}u_{\e})\ran=
g_{\e}^{2-p}|v|^{p-2}{\bf c}\Big[
(1-p)\left(1-2/p\right)g_{\e}^{-2}|v|^{3}\nabla|v|
+\cr
+(p-2)|v|\nabla|v|+v\nabla\overline{v}\Big]\cr
a\lan u_{\e},|u_{\e}|^{p-2}u_{\e}\ran =
ag_{\e}^{2-p}|v|^{p}
}
$$
from which follows
\begin{eqnarray}
\lim_{\e\to 0}\int_{\O}\lan {\bf b}\nabla
u_{\e},|u_{\e}|^{p-2}u_{\e}\ran  =
\int_{\O}\left( -(1-2/p)\, |v|\,
{\bf b}\nabla |v| + \overline{v}\, {\bf b}\nabla v\right)
\label{zac2}\\
\lim_{\e\to 0}\int_{\O}
\lan u_{\e},\overline{\bf c}\nabla(|u_{\e}|^{p-2}u_{\e})\ran  =
\int_{\O}\left( (1-2/p)\, |v| {\bf c}\nabla|v| + v\,
{\bf c}\nabla\overline{v}\right)  \label{zac3}\\
\lim_{\e\to 0}\int_{\O}a\lan u_{\e},|u_{\e}|^{p-2}u_{\e}\ran   =
\int_{\O}a|v|^{2}
\label{zac4}
\end{eqnarray}

From \reff{zac1}--\reff{zac4} we obtain that
$$
\lim_{\e\to 0}\Re\elle(u_{\e},|u_{\e}|^{p-2}u_{\e})
$$
exists and is equal to the left-hand side of \reff{eq:22}. This shows that
\reff{eq:defdis1} implies \reff{eq:22} and so the necessity is proved
for $p\geq 2$.

Let us assume $1<p<2$. Since \reff{eq:defdis2} can be written
as \reff{eq:M}, replacing $\A$, ${\bf b}$, $\overline{\bf c}$
by $\A^{*}$, $-\overline{\bf c}$, $-{\bf b}$ respectively
in formulas \reff{zac1}--\reff{zac4}
we find that
$$
\lim_{\e\to 0}\Re\elle(|u_{\e}|^{p'-2}u_{\e},u_{\e})
$$
exists and is equal to the left-hand side of \reff{eq:P}. Thus
\reff{eq:defdis2} implies \reff{eq:22}.
\end{proof}

\begin{corollary}\label{cor:-1}
    If the form $\elle$ is $L^{p}$-dissipative, we have
   \begin{equation}
        \lan \Re \A \xi,\xi\ran \geq 0
       \label{eq:mainpart}
   \end{equation}
for any $\xi\in\R^{n}$.
\end{corollary}

\begin{proof}
    Given a function $v$,
		   let us set
		     $$
		     X=\Re (|v|^{-1}\overline{v}\, \nabla v), \quad
		     Y= \Im (|v|^{-1}\overline{v}\, \nabla v),
		     $$
		     on the set $\{x\in\O\ |\ v\neq 0\}$. We have
		     $$\display{
		     \Re \lan \A\nabla v,\nabla v\ran=
		     \Re \left\lan \A (|v|^{-1}\overline{v}\, \nabla v),
		     |v|^{-1}\overline{v}\, \nabla v\right\ran=\cr
		     \lan \Re \A X,X\ran +\lan \Re \A Y,Y\ran
		     +\lan \Im(\A-\A^{t})X,Y\ran ,\cr
		     \Re\lan (\A-\A^{*})\nabla(|v|),\nabla v\ran |v|^{-1}v=
		     \Re \lan (\A-\A^{*})X,X+iY\ran = \cr
		     \lan \Im(\A-\A^{*})X,Y\ran ,
		     \cr
		     \Re \lan \A\nabla|v|,\nabla|v|\ran = \lan \Re \A X,X\ran.
		     }
		     $$

   Since $\elle$ is $L^{p}$-dissipative, \reff{eq:22} holds. Hence,
   \beginrighe{uffa2}
   \grande
   \int_{\O}\Big\{ {4\over p\,p'}
   \lan \Re \A X,X\ran + \lan \Re\A Y,Y\ran
	    +\cr
	    \grande 2 \lan (p^{-1}\Im\A + p'^{-1}\Im\A^{*}) X,Y\ran
	    +\lan \Im({\bf b}+{\bf c}),Y\ran |v| +\cr
	    \grande
	    \Re\left[\dive\left( {\bf b}/p - {\bf c}/p'\right) - a
	    \right]|v|^{2}\Big\}   \geq 0
   \endrighe

   We define the
	    function
	    $$
	    v(x)=\r(x)\, e^{i\f(x)}
	    $$
	    where $\r$ and $\f$ are real functions with
	    $\r\in \Cspt^{1}(\O)$ and $\f\in C^{1}(\O)$. Since
	    $$
	   |v|^{-1}\overline{v}\, \nabla v  = |\r|^{-1}(
	    \r \, e^{-i\f}\, (\nabla\r + i\r\nabla\f)\,
	    e^{i\f}) = |\r|^{-1}\r\nabla\r + i|\r|\nabla\f
	    $$
	    on the set $\{x\in\O\ |\ \r(x)\neq 0\}$,
	     it follows from \reff{uffa2} that
	    \beginrighe{3}
		{4\over p\, p'}\,
		\int_{\O}\lan \Re\A \nabla\r,\nabla\r\ran   +
		\int_{\O}\r^{2}\lan \Re\A\nabla\f,\nabla\f\ran    +
		\crg
		2\int_{\O}
		 \r\lan (p^{-1}\Im\A + p'^{-1}\Im\A^{*}) \nabla\r,
\nabla\f\ran
		 +\crg
		 \int_{\O}\r^{}\lan \Im({\bf b}+{\bf c}), \nabla\f \ran+
		 \int_{\O} \Re\left[\dive\left( {\bf b}/p - {\bf
c}/p'\right) - a
	    \right]\r^{2}
		    \geq 0
\endrighe
for any $\r\in \Cspt^{1}(\O)$, $\f\in C^{1}(\O)$.

	    We choose $\f$ by the equality
	    $$
	    \f={\m\over 2}\, \log (\r^{2}+\e)
	    $$
  where $\m\in \R$ and $\e >0$. Then \reff{3} takes the form
	    \beginrighe{3bis}
	    {4\over p\, p'}\,
		    \int_{\O}\lan \Re\A \nabla\r,\nabla\r\ran   +
		    \m^{2}\int_{\O}{\r^{4}\over (\r^{2}+\e)^{2}}\lan
\Re\A\nabla\r,\nabla\r\ran    +
		    \crg
		    2\m \int_{\O}
		     {\r^{2}\over \r^{2}+\e}\lan (p^{-1}\Im\A +
p'^{-1}\Im\A^{*})
		     \nabla\r, \nabla\r\ran
		     +\crg
		     \m\int_{\O}{\r^{3}\over \r^{2}+\e}\lan \Im({\bf
b}+{\bf c}), \nabla\r \ran+
		     \int_{\O} \Re\left[\dive\left( {\bf b}/p - {\bf
c}/p'\right) - a
		\right]\r^{2}
			\geq 0
		     \endrighe

       Letting $\e\to 0^{+}$ in \reff{3bis}
	    leads to
	    \beginrighe{ropositive}
	    {4\over p\, p'}\,
		\int_{\O}\lan \Re\A \nabla\r,\nabla\r\ran   +
		\m^{2}\int_{\O}\lan \Re\A\nabla\r,\nabla\r\ran    +
		\crg
		2\m \int_{\O}
	 \lan (p^{-1}\Im\A + p'^{-1}\Im\A^{*})
		 \nabla\r, \nabla\r\ran
		 +\crg
	 \m\int_{\O}\r\lan \Im({\bf b}+{\bf c}), \nabla\r \ran+
		 \int_{\O} \Re\left[\dive\left( {\bf b}/p - {\bf
c}/p'\right) - a
		    \right]\r^{2}
			    \geq 0.
		 \endrighe

	Since this holds for any $\m\in\R$, we have
\begin{equation}
        \int_{\O}\lan \Re\A\nabla\r,\nabla\r\ran \geq 0
    \label{eq:foranyr}
\end{equation}
	for any $\r\in \Cspt^{1}(\O)$.

	Taking
	$\r(x)=\psi(x)\cos\lan\xi,x\ran$ with a real $\psi\in\Cspt^{1}(\O)$ and
	$\xi\in\R^{n}$, we find
	$$\display{
	\int_{\O}\{ \lan \Re\A\nabla \psi,\nabla\psi\ran \cos^{2}\lan\xi,x\ran -
	[\lan \Re\A \xi,\nabla\psi\ran +\cr
	\lan \Re\A\nabla\psi,\xi\ran ]\, \sin\lan\xi,x\ran\cos\lan\xi,x\ran +
	\lan \Re\A \xi,\xi\ran \psi^{2}(x)\sin^{2}\lan\xi,x\ran\}  \geq 0.
	}
	$$

	On the other hand, taking $\r(x)=\psi(x)\sin\lan\xi,x\ran$,
	$$\display{
	\int_{\O}\{ \lan \Re\A\nabla \psi,\nabla\psi\ran \sin^{2}\lan\xi,x\ran +
	[\lan \Re\A \xi,\nabla\psi\ran + \cr
	\lan \Re\A\nabla\psi,\xi\ran ]\, \sin\lan\xi,x\ran\cos\lan\xi,x\ran
	+\lan \Re\A \xi,\xi\ran \psi^{2}(x)\cos^{2}\lan\xi,x\ran\}  \geq 0.
	}
	$$

	The two inequalities we have obtained lead to
	$$
	\int_{\O}\lan \Re\A\nabla \psi,\nabla\psi\ran   + \int_{\O}\lan \Re\A
	\xi,\xi\ran \psi^{2}  \geq 0.
	$$

	Because of the arbitrariness of $\xi$, we find
	$$
	\int_{\O}\lan \Re\A \xi,\xi\ran \psi^{2}  \geq 0.
	$$

	On the other hand, any nonnegative function $v\in \Cspt(\O)$ can
	be approximated in the uniform norm in $\O$
	by a sequence $\psi_{n}^{2}$, with $\psi_{n}\in
\Cspt^{\infty}(\O)$, and
	then $\lan \Re\A \xi,\xi\ran$ is a nonnegative measure.
\end{proof}

\begin{corollary}\label{cor:0}
    If the form $\elle$ is both $L^{p}$- and $L^{p'}$-dissipative, it is
    also $L^{r}$-dissipa\-tive for any $r$ between $p$ and $p'$, i.e. for
    any $r$ given by
 \begin{equation}
        1/r = t/p + (1-t)/p'
         \qquad (0\leq t\leq 1).
     \label{eq:rinter}
 \end{equation}
\end{corollary}

\begin{proof}
    From the proof of Corollary \ref{cor:-1} we know that
    \reff{uffa2} holds.
In the same way, we find
\beginrighe{uffa3}
\grande
\int_{\O}\Big\{ {4\over p'\,p}
\lan \Re \A X,X\ran + \lan \Re\A Y,Y\ran
	 -\cr
	 \grande 2 \lan (p'^{-1}\Im\A + p^{-1}\Im\A^{*}) X,Y\ran
	 +\lan \Im({\bf b}+{\bf c}),Y\ran |v| +\cr
	 \grande
	 \Re\left[\dive\left(  {\bf b}/p' -  {\bf
c}/p\right) - a
	 \right]|v|^{2}\Big\}   \geq 0.
\endrighe

We multiply \reff{uffa2} by $t$, \reff{uffa3} by $(1-t)$  and
sum up. Since
$$
t/p' + (1-t)/p = 1/r' \quad \hbox{\rm and}\quad
r\, r' \leq p\,p'\, ,
$$
we find, keeping in mind Corollary \ref{cor:-1},
$$\display{
\int_{\O}\Big\{ {4\over r\,r'}
\lan \Re \A X,X\ran + \lan \Re\A Y,Y\ran
	 -\cr
	 \grande 2 \lan (r^{-1}\Im\A + r'^{-1}\Im\A^{*}) X,Y\ran
	 +\lan \Im({\bf b}+{\bf c}),Y\ran |v| +\cr
	+ \Re\left[\dive\left( {\bf b}/r -  {\bf
c}/r'\right) - a
	 \right]|v|^{2}\Big\}   \geq 0
}
$$
and $\elle$ is $L^{r}$-dissipative by Lemma \ref{lemma:1} .
\end{proof}

\begin{corollary}
     Suppose that either
  \begin{equation}
        \Im\A=0, \qquad \Re\dive
          {\bf b}= \Re\dive{\bf c}=0
      \label{eq:uffa}
  \end{equation}
  or
  \begin{equation}
      \Im\A=\Im\A^{t},\quad \Im({\bf b} +{\bf c})=0, \quad \Re\dive
          {\bf b}= \Re\dive{\bf c}=0.
       \label{eq:uffa2}
   \end{equation}
    If $\elle$ is $L^{p}$-dissipative, it is also $L^{r}$-dissipative
    for any $r$ given by \reff{eq:rinter}.
\end{corollary}

\begin{proof}
	     Assume that \reff{eq:uffa} holds.
	     With the notation introduced in Corollary \ref{cor:-1},
	 inequality \reff{eq:22} reads as
	 $$\display{
	 \int_{\O}\Big({4\over p\,p'}\lan \Re \A X,X\ran + \lan \Re\A Y,Y\ran
	 +\cr
	 \lan \Im({\bf b}+{\bf c}),Y\ran |v| - \Re a |v|^{2}\Big)
	   \geq 0.
	 }
	 $$
  Since the left-hand side
	does not change after replacing $p$ by $p'$,  Lemma \ref{lemma:1}
	gives the result.

 Let \reff{eq:uffa2} holds. Using the formula
	 \beginrighe{giochetto}
	        p^{-1}\Im\A +p'^{-1}\Im\A^{*}=\crg
		p^{-1}\Im\A
	             -p'^{-1}\Im\A^{t}=-(1-2/p)\Im\A,
	 \endrighe
	    we obtain
    $$\display{
    \int_{\O}\Big( {4\over p\,p'}
    \lan \Re \A x,x\ran + \lan \Re\A Y,Y\ran
	     -\cr
	     2\, (1-2/p) \lan \Im\A X,Y\ran
	     - \Re a
	    |v|^{2}\Big)   \geq 0.
	     }
	     $$

	     Replacing $v$ by $\overline{v}$, we find
	     $$\display{
		 \int_{\O}\Big( {4\over p\,p'}
		 \lan \Re \A x,x\ran + \lan \Re\A Y,Y\ran
			  +\cr
			  \grande 2\, (1-2/p) \lan \Im\A X,Y\ran
			  \grande
			  - \Re a
			 |v|^{2}\Big)   \geq 0
			  }
			  $$
  and we have the $L^{p'}$-dissipativity by $1-2/p = -1+2/p'$. 
   The reference to Corollary
   \ref{cor:0} completes the proof.
\end{proof}

We give now a sufficient condition for the $L^{p}$-dissipativity.
This is a direct consequence of Lemma \ref{lemma:1}.

\begin{corollary}\label{cor:1}
    Let $\a,\b$ two real constants. If
   \beginrighe{polyn}
    {4\over p\,p'}\lan \Re \A\xi,\xi\ran + \lan \Re \A\h,\h\ran
       +2 \lan (p^{-1}\Im\A + p'^{-1}\Im\A^{*}) \xi,\h \ran +\crg
       \lan \Im({\bf b} + {\bf c}),\h\ran 
       -
       2 \lan \Re (\a{\bf b}/p - \b {\bf c}/p'),\xi\ran
       +\cr
       \Re\left[\dive\left((1-\a){\bf b}/p -  (1-\b){\bf c}/p'\right) - a
       \right]\geq 0
\endrighe
       for any $\xi,\h\in\R^{n}$,
    the form $\elle$ is $L^{p}$-dissipative.
\end{corollary}

\begin{proof}
    In the proof of Lemma \ref{lemma:1} we have integrated by parts 
	in \reff{eq:part1} and \reff{eq:part2}. More generally, 
	we have
	$$\display{
	2/p\int_{\O}\lan \Re{\bf b},\Re(\overline{v}\nabla v)\ran =
	2\a/p \int_{\O}\lan \Re{\bf b},\Re(\overline{v}\nabla v)\ran
	-\cr
	(1-\a)/p\int_{\O}\Re(\nabla^{t}{\bf b})|v|^{2}\, ;\cr
	2/p'\int_{\O}\lan \Re{\bf c}, \Re(\overline{v}\nabla v)\ran =
	2\b/p' \int_{\O}\lan \Re{\bf c}, \Re(\overline{v}\nabla 
	v)\ran-\cr
	(1-\b)/p' \int_{\O}\Re(\nabla^{t}{\bf c})|v|^{2}\, .
	}
	$$
	
	This leads to write conditions \reff{eq:22} in a slightly different 
	form:
	$$\display{
		\Re \int_{\O}\Big[ \lan \A\nabla v,\nabla v\ran -
		(1-2/p)\lan (\A-\A^{*})\nabla(|v|),|v|^{-1}\overline{v}\nabla v\ran  -
		\cr
		(1-2/p)^{2}\lan \A \nabla(|v|),\nabla(|v|)\ran
	\Big]
	+
	\int_{\O}\lan\Im ({\bf b}+{\bf c}),\Im(\overline{v}\nabla
	v)\ran\,   -
	\cr
	2\int_{\O}\lan  \Re (\a{\bf b}/p - \b {\bf c}/p'), 
	\Re(\overline{v}\nabla v) \ran +\cr
	\int_{\O}\Re(\dive ((1-\a){\bf b}/p - (1-\b){\bf c}/p') - a
	)|v|^{2}
		\geq 0.
		}
	$$
    
By using the functions $X$ and $Y$ introduced in Corollary
	  \ref{cor:-1},
	 the left-hand side of the last inequality can be written as
	  $$
	  \int_{\O} Q(X,Y)
	  $$
	  where $Q$ denotes the polynomial \reff{polyn}.
	  The result follows  from Lemma \ref{lemma:1}.
    \end{proof}

    Generally speaking,  conditions \reff{polyn} are not
    necessary for $L^{p}$-dissipa\-tivity. We show this by the following 
    example,
    where
    $\Im\A$ is not symmetric. Later we 
    give another example showing that, even for symmetric matrices
    $\Im\A$,
    conditions \reff{polyn} are not necessary 
    for $L^{p}$-dissipa\-tivity (see Example \ref{ex:1}).
Nevertheless in the next section
    we show that the conditions are necessary for the $L^{p}$-dissipativity, 
    provided the operator
    $A$ has no lower order terms and the matrix
    $\Im\A$ is symmetric (see Theorem \ref{th:1new} and Remark \ref{rm:1}).

    \begin{example}
        Let $n=2$ and
	$$
	\A=\left(\begin{array}{cc}
	1 & i\g \\ -i\g & 1
	\end{array}\right)
	$$
	where $\g$ is a real constant, ${\bf b}={\bf c}=a=0$. In this case
polynomial \reff{polyn}
	is given by
	$$
	(\h_{1}-\g\xi_{2})^{2}+	(\h_{2}-\g\xi_{1})^{2}
-(\g^{2}-4/(pp'))|\xi|^{2}.
	$$
	
	Taking $\g^{2} >
	4/(pp')$,  condition \reff{polyn} is not satisfied,
	while we have the $L^{p}$-dissipativity, because the corresponding
operator $A$
	is the Laplacian.
    \end{example}

\paragrafo{The operator $\dive(\A\nabla u)$}\label{sec:main}

In this section we  consider operator \reff{eq:N} without lower order
terms:
\begin{equation}
    Au=\dive(\A\nabla u)
    \label{eq:nolower}
\end{equation}
with the coefficients $a^{hk}\in
(\Cspt(\O))^{*}$.
The following Theorem  contains 
an algebraic necessary and sufficient condition for the
$L^{p}$-dissipativity. 

This result is new even for
smooth coefficients, when it implies a criterion for the
$L^{p}$-contractivity of the corresponding semigroup (see Theorem
\ref{th:contr} below).

\begin{theorem}\label{th:1new}
    Let the matrix $\Im\A$ be symmetric, i.e. $\Im\A^{t}=\Im\A$. The  form
$$
        \elle(u,v)=\int_{\O}\lan \A\nabla u,\nabla v\ran
$$
    is $L^{p}$-dissipative if and only if
    \begin{equation}\label{eq:24}
	|p-2|\, |\lan \Im\A\xi,\xi\ran| \leq 2 \sqrt{p-1}\,
	\lan \Re\A\xi,\xi\ran
    \end{equation}
	for any $\xi\in\R^{n}$, where $|\cdot|$ denotes the total
	variation.
\end{theorem}

\begin{proof}

    {\it Sufficiency.} In view of Corollary
    \ref{cor:1}
    the form $\elle$ is $L^{p}$-dissipative if
\begin{equation}
       {4\over p\,p'}\lan\Re\A\xi,\xi\ran + \lan \Re\A\h,\h\ran
       -2(1-2/p)\lan\Im\A\xi,\h\ran \geq 0
    \label{eq:25}
\end{equation}
    for any $\xi,\h\in\R^{n}$.

    By putting
	      $$
	      \lambda = {2\sqrt{p-1}\over p}\ \xi
	      $$
	      we write \reff{eq:25} in the form
	      $$
	      \lan \Re \A\lambda,\lambda\ran +
	      \lan \Re \A\h,\h\ran
		      -{p-2\over \sqrt{p-1}}\lan
		      \Im\A\lambda,\h)\geq 0.
	      $$

	      Then \reff{eq:25} is equivalent to
	      $$
		  {\mathscr S}(\xi,\h) :=
		  \lan \Re \A\xi,\xi\ran +
			    \lan \Re \A\h,\h\ran -
				    {p-2\over \sqrt{p-1}}\lan
				    \Im \A \xi,\h)
		      \geq 0
		      $$
    for any $\xi,\h\in\R^{n}$.

	      For any nonnegative $\f\in\Cspt(\O)$, define
	      $$
	      \l_{\f} = \min_{|\xi|^{2}+|\h|^{2}=1}\int_{\O}{\mathscr S}
(\xi,\h)
	      \, \f\,  .
	      $$

	      Let us fix $\xi_{0},\h_{0}$ such that
$|\xi_{0}|^{2}+|\h_{0}|^{2}=1$ and
	      $$
	     \l_{\f}=\int_{\O}{\mathscr S} (\xi_{0},\h_{0})
			\, \f\, .
			$$
		     We have the algebraic system
	      $$
	      \cases{\grande
	      \int_{\O}\left(2\, \Re\A\xi_{0}-{p-2\over 2\, \sqrt{p-1}}\,
	      \Im(\A-\A^{*})\h_{0}\right)\f\,   =2\,
	      \l_{\f}\, \xi_{0}\cr
	      \grande
	      \int_{\O}\left(2\, \Re\A\h_{0}-{p-2\over 2\, \sqrt{p-1}}\,
			\Im(\A-\A^{*})\xi_{0}\right)\f\,   =2\,
			\l_{\f}\, \h_{0}\, .
	      }
	      $$

	      This implies
	      $$\display{
	      \int_{\O}\left(2\, \Re\A(\xi_{0}-\h_{0}) + {p-2\over 2\,
\sqrt{p-1}}\,
			\Im(\A-\A^{*})(\xi_{0}-\h_{0})\right)\f\,   =\cr
			2\,
			\l_{\f}\, (\xi_{0}-\h_{0})
			}
	      $$
	      and therefore
	      $$\display{
	     \int_{\O}\left( 2\lan \Re\A(\xi_{0}-\h_{0}),\xi_{0}-\h_{0}\ran
+ {p-2\over \sqrt{p-1}}\lan
				  \Im\A(\xi_{0}-\h_{0}),
\xi_{0}-\h_{0}\ran\right)\f\,   = \cr
				  2\, \l_{\f}\, |\xi_{0}-\h_{0}|^{2}.
				  }
	      $$

	      The left-hand side is nonnegative because of \reff{eq:24}. Hence,
	      if $\l_{\f}<0$, we find $\xi_{0}=\h_{0}$. On the
	      other hand we have
	      $$\display{
	      \l_{\f}=\int_{\O}{\mathscr S}(\xi_{0},\xi_{0})\, \f\,   = \cr
	      \int_{\O}\left(2 \lan \Re\A\xi_{0},\xi_{0}\ran  -
		      {p-2\over \sqrt{p-1}} \lan
\Im\A\xi_{0},\xi_{0}\ran\right) \f\,
		      \geq 0.
		      }
	      $$

	      This shows that  $\l_{\f}\geq 0$ for any nonnegative $\f$ and
the
	      sufficiency is proved.

    {\it Necessity.}
    We know from the proof of  Corollary \ref{cor:-1} that if $\elle$
    is $L^{p}$-dissipative, then \reff{ropositive} holds for any
    $\r\in\Cspt^{1}(\O)$, $\m\in\R$. In the
    present case, keeping in mind \reff{giochetto}, \reff{ropositive} can
be written as
    $$
    \int_{\O}\lan \B\nabla\r,\nabla\r\ran   \geq 0,
    $$
     where
$$
\B={4\over p\,p'}\Re\A + \m^{2}\Re\A-2\, \m\, (1-2/p)\Im\A
$$

In the proof of Corollary
\ref{cor:-1}, we have also seen  that from \reff{eq:foranyr} for any
$\r\in\Cspt^{1}(\O)$,
\reff{eq:mainpart} follows. In the same way,  the last relation
implies
$\lan\B\xi,\xi\ran \geq 0$, i.e.
	  $$\display{
	  {4\over p\, p'}\,
		  \lan \Re\A \xi,\xi\ran  +
		  \m^{2} \lan \Re\A \xi,\xi\ran -
		 2\, \m\, (1-2/p)\lan \Im \A \xi,\xi\ran
		  \geq 0 }
		  $$
	  for any $\xi\in \R^{n}$, $\m\in \R$.

	  Because of the arbitrariness of $\m$ we have
	 $$\display{
	 \int_{\O}\lan \Re\A \xi,\xi\ran\, \f\,   \geq 0
	 \cr
	 (1-2/p)^{2}\left(\int_{\O} \lan \Im \A
	   \xi,\xi\ran \,
	   \f\,  \right)^{2}
	  \leq {4\over p\, p'} \left(\int_{\O}\lan \Re\A \xi,\xi\ran \,\f\,
	   \right)^{2}
	  }
	  $$
	  i.e.
$$
	       |p-2| \  \left|\int_{\O} \lan \Im \A \xi,\xi\ran \,\f\,
\right| \leq
	      2\  \sqrt{p-1}\ \int_{\O}\lan \Re \A \xi,\xi\ran \,\f\,
$$
for any $\xi\in\R^{n}$ and for any nonnegative $\f\in\Cspt(\O)$.

We have
   $$
   |p-2|\, \left|\int_{\O}\lan \Im\A\xi,\xi\ran \f 
\right| \leq 2 \sqrt{p-1}\, \int_{\O}\lan \Re\A\xi,\xi\ran |\f| 
   $$
   for any $\f\in \Cspt(\O)$ and this implies \reff{eq:24}, because
$$\display{
|p-2| \int_{\O}|\lan \Im\A\xi,\xi\ran|\, g =|p-2| \sup_{\f\in \Cspt(\O)\atop |\f|\leq g}
   \left| \int_{\O}\lan \Im\A\xi,\xi\ran 
  \f \right| \leq \cr
  2\sqrt{p-1} 
  \sup_{\f\in \Cspt(\O)\atop |\f|\leq g}\int_{\O}\lan \Re\A\xi,\xi\ran |\f|
  \leq 
  2\sqrt{p-1} \int_{\O}\lan \Re\A\xi,\xi\ran g 
  }
$$
for any nonnegative $g\in \Cspt(\O)$.
\end{proof}

\begin{remark} \label{rm:1} From the proof of  Theorem \ref{th:1new} we
see that condition \reff{eq:24} holds if and only if
$$
{4\over p\,p'}\lan \Re\A\xi,\xi\ran +\lan \Re\A\h,\h\ran
-2(1-2/p)\lan\Im\A\xi,\h\ran \geq 0
$$
for any $\xi,\h\in\R^{n}$. This means that   conditions
\reff{polyn} are
necessary and sufficient for the operators
considered in Theorem \ref{th:1new}.
\end{remark}

\begin{remark}
    Let us assume that either $A$ has lower order terms
	or they are absent and $\Im\A$ is not symmetric.
Using the same arguments as in Theorem \ref{th:1new},
one could prove that \reff{eq:24} is
still a necessary condition
for $A$ to be $L^{p}$-dissipative. However, in general, it is not
sufficient. This is shown by the next example (see also
Theorem \ref{th:const} below  for the particular case of constant
coefficients).
\end{remark}

\begin{example} Let $n=2$ and let $\O$ be a bounded domain. Denote by $\s$ 
a  not identically vanishing real function in $\Cspt^{2}(\O)$ and let 
$\l\in\R$.
Consider operator \reff{eq:nolower} with
$$
\A=\left(\begin{array}{cc}
	1 & i\l\de_{1}(\s^{2}) \\ -i\l\de_{1}(\s^{2}) & 1
	\end{array}\right)
$$
i.e.
$$
Au=\de_{1}(\de_{1}u+i\l\de_{1}(\s^{2})\, \de_{2}u) + 
\de_{2}(-i\l\de_{1}(\s^{2})\,
\de_{1}u+\de_{2}u),
$$
where $\de_{i}=\de/\de x_{i}$ ($i=1,2$).

By definition, we have $L^{2}$-dissipativity if and only if
$$
\Re\int_{\O}((\de_{1}u+i\l\de_{1}(\s^{2})\, \de_{2}u)
\de_{1}\overline{u}
+ (-i\l\de_{1}(\s^{2})\,
\de_{1}u+\de_{2}u)\de_{2}\overline{u})\, dx \geq 0
$$
for any $u\in\Cspt^{1}(\O)$, i.e. if and only if
$$
\int_{\O}|\nabla u|^{2}dx -2\l\int_{\O}\de_{1}(\s^{2})
\Im(\de_{1}\overline{u}\,\de_{2}u)\, dx \geq 0
$$
for any $u\in\Cspt^{1}(\O)$.
Taking $u=\s\, \exp(itx_{2})$ ($t\in\R$), we obtain, in particular,
\begin{equation}
    t^{2}\int_{\O}\s^{2}dx -t\l 
    \int_{\O}(\de_{1}(\s^{2}))^{2} dx  + \int_{\O}|\nabla \s|^{2}dx \geq 0.
    \label{eq:inpartic}
\end{equation}

Since
$$
\int_{\O}(\de_{1}(\s^{2}))^{2} dx > 0,
$$
we can choose $\l\in\R$ so that
\reff{eq:inpartic} is impossible
for all $t\in\R$.
Thus  $A$ is not $L^{2}$-dissipative, although
\reff{eq:24} is satisfied.

Since $A$ can be written as
$$
	Au=\D u - i\l(\de_{21}(\s^{2})\,
	\de_{1}u - \de_{11}(\s^{2})\, \de_{2}u),
$$
the same example shows
that \reff{eq:24} is not sufficient for the $L^{2}$-dissipativity
in the presence of lower order terms, even if $\Im\A$ is symmetric.
\end{example}

\paragrafo{General equation with constant 
coefficients}\label{sec:constant}

In this section we  characterize the $L^{p}$-dissipativity for
a differential operator $A$, say
\begin{equation}
    Au=\nabla^{t}(\A\nabla u) + {\bf b}\nabla u + au
    \label{eq:genform}
\end{equation}
with constant complex coefficients. Without loss of generality
 we 
assume that the matrix $\A$ is symmetric.

\begin{theorem}\label{th:const}
    Let $\O$ be an open set in $\R^{n}$ which contains balls of
    arbitrarily large radius. The operator $A$
    is $L^{p}$-dissipative if and only if there exists a real constant
    vector $V$ such that
    \begin{eqnarray}
        2\Re \A V+\Im {\bf b}=0
        \label{eq:sistV}\\
	\Re a +  \lan \Re\A V,V\ran \leq 0 &
	    \label{eq:V2}
    \end{eqnarray}
    and the inequality
    \begin{equation}\label{eq:V1}
	|p-2|\, |\lan \Im\A\xi,\xi\ran| \leq 2 \sqrt{p-1}\,
	\lan \Re\A\xi,\xi\ran
    \end{equation}
   holds for any $\xi\in\R^{n}$.
\end{theorem}

\begin{proof}
    First, let us prove the Theorem for the special case ${\bf
    b}=0$, i.e. for the operator
    $$
    A=\nabla^{t}(\A\nabla u)  + au.
    $$

  If  $A$ is $L^{p}$-dissipative, \reff{eq:22} holds  for
  any $v\in \Cspt^{1}(\O)$.  We find, by  repeating 
  the arguments used in
  the proof of Theorem \ref{th:1new}, that
 \beginrighe{last}
  {4\over p\,p'}\int_{\O}\lan \Re\A\nabla \r, \nabla \r \ran\, dx +
  \m^{2}\int_{\O}\lan \Re\A\nabla \r, \nabla \r \ran\, dx -\crg
  2\, \m\, (1-2/p)\int_{\O}\lan \Im\A\nabla \r, \nabla \r \ran\, dx
  - (\Re a)\int_{\O}\r^{2}dx \geq 0
\endrighe
 for any $\r\in\Cspt^{\infty}(\O)$ and for any $\m\in\R$. As in the proof
of Theorem
 \ref{th:1new} this implies \reff{eq:V1}. On the other hand, 
we can find
 a sequence of balls contained in
 $\O$ with centres $x_{m}$ and radii $m$.  Set
 $$
\r_{m}(x)=m^{-n/2}\s\left((x-x_{m})/ m\right),
$$
 where $\s\in\Cspt^{\infty}(\R^{n})$, $\spt\s\subset B_{1}(0)$ and
 $$
 \int_{B_{1}(0)}\s^{2}(x)\, dx =1.
 $$
 
 Putting in \reff{last} $\m=1$ and $\r=\r_{m}$, we obtain
$$\display{
  {4\over p\,p'}\int_{B_{1}(0)}\lan \Re\A\nabla \s, \nabla \s \ran\, dy +
 \int_{B_{1}(0)}\lan \Re\A\nabla \s, \nabla \s \ran\, dy \ -\cr
  2\,  (1-2/p)\int_{B_{1}(0)}\lan \Im\A\nabla \s, \nabla \s \ran\, dy
  - m^{2}(\Re a) \geq 0
  }
  $$
for any $m\in\N$. This implies $\Re a \leq 0$. Note that in this
case the algebraic system \reff{eq:sistV}  has always the trivial
solution and that for any eigensolution $V$ (if they  exist)
we  have $\lan \Re\A V,V\ran=0$. Then \reff{eq:V2} is satisfied.

  Conversely, if \reff{eq:V1} is satisfied, we have (see Remark
  \ref{rm:1})
  $$
        {4\over p\,p'}\lan \Re\A\xi,\xi\ran + \lan \Re \A\h,\h\ran
        - 2\, (1-2/p) \lan \Im \A\xi,\xi\ran \geq 0
$$
   for any $\xi,\h\in\R^{n}$. If also \reff{eq:V2} is satisfied
   (i.e. if $\Re a \leq 0$), $A$ is $L^{p}$-dissipative  in view of Corollary
\ref{cor:1}.

Let us consider the operator in the general form \reff{eq:genform}.
If $A$ is $L^{p}$-dissipative, we find, by repeating the arguments employed in
the proof of Theorem \ref{th:1new}, that 
$$\display{
  {4\over p\,p'}\int_{\O}\lan \Re\A\nabla \r, \nabla \r \ran\, dx +
  \int_{\O}\r^{2}\lan \Re\A\nabla \f, \nabla \f \ran\, dx -\cr
 2\, (1-2/p)\int_{\O}\r\,\lan \Im\A\nabla \r, \nabla \f \ran\, dx
  +\cr
  \int_{\O}\r^{2}\lan \Im {\bf b},\nabla\f\ran\, dx - \Re
a\int_{\O}\r^{2}dx \geq 0
  }
  $$
for any $\r\in \Cspt^{1}(\O)$, $\f\in C^{1}(\O)$. By
fixing $\r$ and choosing $\f=t\lan \h, x\ran$ ($t\in \R$, $\h\in\R^{n}$) we get
$$\display{
  {4\over p\,p'}\int_{\O}\lan \Re\A\nabla \r, \nabla \r \ran\, dx +
 ( t^{2}\lan \Re \A \h,\h\ran +t\, \lan \Im {\bf b},\h\ran - \Re a)
 \int_{\O}\r^{2}\, dx  \geq 0
  }
  $$
for any $t\in\R$. This leads to
$$
|\lan \Im {\bf b},\h\ran |^{2}\leq K\,\lan\Re\A\h,\h\ran
$$
for any $\h\in\R^{n}$ and this inequality shows that system \reff{eq:sistV} is
solvable. Let $V$ be a solution of this system and  let
$$
z=e^{-i\lan V,x\ran}u.
$$

One checks directly that
$$
Au=(\nabla^{t}(\A\nabla z) + \lan {\bf c},\nabla z\ran +
\a z)e^{i\lan V,x\ran}
$$
where
$$
{\bf c}=2i \A V+{\bf b}, \quad
\a=a+i\lan {\bf b},V\ran - \lan \A V,V\ran.
$$

Since we have
$$
\int_{\O}\lan Au,u\ran|u|^{p-2}dx=\int_{\O}\lan
\nabla^{t}(\A\nabla z) + \lan {\bf c},\nabla z\ran +
\a z,z\ran |z|^{p-2}dx\, ,
$$
the $L^{p}$-dissipativity of $A$ is equivalent to the $L^{p}$-dissipativity of the
operator
$$
\nabla^{t}(\A\nabla z) + \lan {\bf c},\nabla z\ran +
\a z\, .
$$

On the other hand Lemma \ref{lemma:1} shows that, as far as the first
order terms are concerned,  the $\Re{\bf b}$ does not play any role.
Since
$
\Im \bf c=0
$
because of \reff{eq:sistV}, the $L^{p}$-dissipativity of $A$ is
equivalent to the $L^{p}$-dissipativity of the operator
\begin{equation}
    \nabla^{t}(\A\nabla z) +
\a z\,.
    \label{eq:opeq}
\end{equation}

By what we have already proved above,  the last
operator is $L^{p}$-dissipative if and only if \reff{eq:V1} is satisfied and
$\Re\a\leq 0$. From \reff{eq:sistV} it follows
that $\Re\a$ is equal to the left-hand side of \reff{eq:V2}.

Conversely, if there exists a solution $V$ of \reff{eq:sistV},
\reff{eq:V2}, and if \reff{eq:V1} is satisfied,  operator
\reff{eq:opeq} is $L^{p}$-dissipative. Since this is equivalent to the
$L^{p}$-dissipativity of $A$, the proof is complete.
\end{proof}

\begin{corollary}\label{cor:const}
    Let $\O$ be an open set in $\R^{n}$ which contains balls of
       arbitrarily large radius. Let us suppose that the matrix $\Re\A$ 
       is
       not singular. The operator $A$
       is $L^{p}$-dissipative if and only if \reff{eq:V1} holds and
       \begin{equation}
           4\Re a \leq -\lan (\Re\A)^{-1}\Im {\bf b}, \Im {\bf b}\ran
           \label{eq:newV2}
       \end{equation}
\end{corollary}

\begin{proof}
    If $\Re\A$ is not singular, the only vector $V$ satisfying
    \reff{eq:sistV} is
    $$
    V=-(1/2) (\Re\A)^{-1}\Im {\bf b}
    $$
    and \reff{eq:V2} is satisfied if and only if
\reff{eq:newV2} holds. The result follows from Theorem
\ref{th:const}.
\end{proof}

\begin{example}\label{ex:1}
Let $n=1$ and $\O=\R^{1}$. Consider the operator
$$
\left( 1 + 2\, {\sqrt{p-1}\over p-2}\, i\right) u'' +2i u' -u,
$$
where $p\neq 2$ is fixed. Conditions
\reff{eq:V1} and \reff{eq:newV2} are satisfied
 and this operator is $L^{p}$-dissipative, in view of Corollary \ref{cor:const}.

On the other hand, the polynomial considered in Corollary \ref{cor:1} is
$$
Q(\xi,\h)=\left(2\, {\sqrt{p-1}\over p}\, \xi -\h\right)^{2}
+2\h+1
$$
which is not nonnegative for any $\xi,\h\in\R$.
This shows that, in general, condition \reff{polyn} is not necessary for the
$L^{p}$-dissipativity, even if the matrix $\Im\A$ is symmetric.
\end{example}

\paragrafo{Smooth coefficients}\label{sec:smooth}

Let us consider the operator
\begin{equation}
    Au = \dive(\A\nabla u) + {\bf b}\nabla u +a\, u
    \label{eq:Acompleto}
\end{equation}
with the coefficients $a^{hk}, b^{h}\in C^{1}(\overline{\O})$, $a\in
C^{0}(\overline{\O})$. Here $\O$ is a bounded domain in $\R^{n}$,
whose boundary is in the class $C^{2,\a}$ for some $\a\in[0,1)$
(this regularity assumption could be weakened, but we prefer to avoid 
the technicalities related to such generalizations).

We consider $A$ as an operator defined on the set
\begin{equation}
    \Dom(A)=W^{2,p}(\O)\cap \Wspt^{1,p}(\O).
    \label{eq:DomA}
\end{equation}

\begin{definition} The operator $A$ is said to be $L^{p}$-dissipative if
\begin{equation}
    \Re \int_{\O}\lan Au,u\ran |u|^{p-2}dx \leq 0
    \label{eq:defclass}
\end{equation}
for any $u\in\Dom(A)$.
\end{definition}

We  show that  the  $L^{p}$-dissipativity of $A$ is
equivalent to the $L^{p}$-dissipativity of the sesquilinear form 
$$
    {\elle}(u,v)=\int_{\O}(\lan \A\nabla u,\nabla v\ran -\lan {\bf b}\nabla u,v\ran
    -a\lan u,v\ran)\,
$$

\begin{lemma}\label{lemma:2}
    The form $\elle$ is $L^{p}$-dissipative if and only if
   \beginrighe{vpos}
        \Re \int_{\O}\Big[ \lan \A\nabla v,\nabla v\ran -
               (1-2/p)\lan (\A-\A^{*})\nabla(|v|),|v|^{-1}\overline{v}\nabla v\ran  -\\
               \grande
               (1-2/p)^{2}\lan \A \nabla(|v|),\nabla(|v|)\ran
           \Big]dx +\crg
	   \int_{\O}\lan \Im {\bf b},\Im (\overline{v}\nabla v)\ran dx +
	   \int_{\O}\Re(\nabla^{t}({\bf b}/p)-a)|v|^{2}dx
	   \geq 0
   \endrighe
     for any $v\in \Hspt^{1}(\O)$.
\end{lemma}

\begin{proof}
    
    \textit{Sufficiency.}
    We know from Lemma \ref{lemma:1} that $\elle$ is
    $L^{p}$-dissipative if and only if \reff{vpos} holds for any
    $v\in\Cspt^{1}(\O)$. Since $\Cspt^{1}(\O)\subset \Hspt^{1}(\O)$, 
  the sufficiency follows.

    \textit{Necessity.}
   Given $v\in \Hspt^{1}(\O)$, we can find a 
   sequence $\{v_{n}\}\subset \Cspt^{1}(\O)$ such that $v_{n}\to v$ in $\Hspt^{1}(\O)$.
   Let us show that
   \begin{equation}
       \chi_{E_{n}}|v_{n}|^{-1}\overline{v}_{n}\nabla v_{n} \to
       \chi_{E}|v|^{-1}\overline{v}\nabla v 
       \quad \hbox{\rm in}\ L^{2}(\O)
       \label{eq:convl2}
   \end{equation}
   where $E_{n}=\{x\in\O\ |\ v_{n}(x)\neq 0\}$, $E=\{x\in\O\ |\ 
   v(x)\neq 0\}$. We may  assume 
   $v_{n}(x)\to v(x)$, $\nabla v_{n}(x) \to \nabla v(x)$ almost 
   everywhere in $\O$. We see that
   \begin{equation}
       \chi_{E_{n}}|v_{n}|^{-1}\overline{v}_{n}\nabla v_{n} \to
	     \chi_{E} |v|^{-1}\overline{v}\nabla v 
       \label{eq:conqo}
   \end{equation}
   almost everywhere on the set $E\cup \{x\in \O\setminus E\ |\ 
   \nabla v(x)=0\}$. Since the set $\{x\in \O\setminus E\ |\ 
   \nabla v(x)\neq 0\}$ has zero measure, we can say that 
   \reff{eq:conqo} holds almost everywhere in $\O$.
   
   Moreover, since
   $$
   \int_{G}|\chi_{E_{n}}|v_{n}|^{-1}\overline{v}_{n}\nabla v_{n} 
   |^{2}dx \leq \int_{G}|\nabla v_{n}|^{2}dx
   $$
   for any measurable set $G\subset \O$ and $\{\nabla v_{n}\}$ is  
   convergent in $L^{2}(\O)$,  the sequence
   $\{|\, \chi_{E_{n}}|v_{n}|^{-1}\overline{v}_{n}\nabla v_{n} -
      \chi_{E}|v|^{-1}\overline{v}\nabla v\,|^{2}\}$  has 
   uniformly absolutely continuos integrals. 
   Now we may appeal to Vitali's Theorem to obtain \reff{eq:convl2}.

   From this it follows that \reff{vpos} for any $v\in \Hspt^{1}(\O)$
   implies \reff{vpos} for any $v\in\Cspt^{1}(\O)$.
   Lemma \ref{lemma:1} shows that $\elle$ is $L^{p}$-dissipative.
   \end{proof}

   \begin{lemma}\label{lemma:equiv}
       The
       form $\elle$ is $L^{p}$-dissipative if and only if
       \begin{equation}
           \Re\int_{\O}(\lan \A \nabla u,\nabla(|u|^{p-2}u)\ran -
	   \lan {\bf b}\nabla u, |u|^{p-2}u\ran -a\, |u|^{p})dx
	   \geq 0
           \label{upos}
       \end{equation}
       for any $u\in \Xi$, where $\Xi$ denotes the space $\{u\in
       C^{2}(\overline{\O})\ |\ u|_{\de\O}=0\}$.
   \end{lemma}

   \begin{proof}

       \textit{Necessity.} Since $\elle$ is $L^{p}$-dissipative,
       \reff{vpos} holds for any $v\in\Hspt^{1}(\O)$.
       Let $u\in \Xi$. We introduce the function
       $$
	\r_{\e}(s)=\cases{\e^{p-2\over 2} & if $0\leq s\leq \e$\cr
	s^{p-2\over 2} & if $s>\e$}
       $$

       Setting
       $$
       v_{\e}=\r_{\e}(|u|)\, u
       $$
       a direct computation shows that $u=\s_{\e}(|v_{\e}|)\, v_{\e}$ 
       and
       $\r^{2}_{\e}(|u|)\, u\- =[\s_{\e}(|v_{\e}|)]^{-1}\, v_{\e}$,
       where
       $$
       \s_{\e}(s)=\cases{\e^{2-p\over 2} & if $0\leq s\leq \e^{p\over 2}$\cr
	s^{2-p\over p} & if $s>\e^{p\over 2}$.}
       $$

       Therefore
       $$\displaylines{
        \langle\A \nabla u, \nabla[\r^{2}_{\e}(|u|)\,
       u]\rangle=
       \langle \A\nabla[\s_{\e}(|v_{\e}|)\, v_{\e}],
       \nabla [ (\s_{\e}(|v_{\e}|))^{-1}{{v_{\e}}}] \rangle =\cr
      \langle \A \left[\s_{\e}(|v_{\e}|)\, \nabla v_{\e} +
       \s_{\e}'(|v_{\e}|)\, v_{\e}\, \nabla |v_{\e}|\right]
       , \s_{\e}(|v_{\e}|)^{-1} \nabla {v_{\e}} -
       \cr
       \s_{\e}'(|v_{\e}|) \s_{\e}^{-2}(|v_{\e}|){v_{\e}}\,
       \nabla |v_{\e}|  \rangle  = \cr
       \langle \A\nabla v_{\e}, \nabla v_{\e}\rangle +
       \s_{\e}'(|v_{\e}|) \s_{\e}(|v_{\e}|)^{-1} \left(
        \langle v_{\e} \A  \nabla |v_{\e}|, \nabla v_{\e}\rangle-
        \langle \A\nabla v_{\e},v_{\e}\, \nabla |v_{\e}|\rangle\right)-\cr
       -\s_{\e}'(|v_{\e}|)^{2} \s_{\e}(|v_{\e}|)^{-2} \langle
       v_{\e}\A\nabla|v_{\e}|, v_{\e}\nabla|v_{\e}|\rangle\, .
       \cr}
       $$

       Since
       $$
       {\s_{\e}'(|v_{\e}|)\over \s_{\e}(|v_{\e}|)}=\cases{
       0 & if $0<|u|<\e$\cr
       -(1-2/p)\ |v_{\e}|^{-1} & if $|u| >\e$}
       $$
       we may write
       $$\displaylines{
       \int_{\O}\langle \A\nabla u, \nabla [\r^{2}_{\e}(|u|)\,
       u]\rangle\, dx =
       \int_{\O} \langle \A\nabla v_{\e}, \nabla v_{\e}\rangle\, dx -\cr
       -(1-2/p) \int_{E_{\e}}{1\over |v_{\e}|} \left(
        \langle v_{\e}\A \nabla |v_{\e}|, \nabla v_{\e}\rangle-
      \langle \A\nabla v_{\e},v_{\e}\, \nabla |v_{\e}|\rangle\right)\, dx -\cr
       - ( 1-2/p)^{2} \int_{E_{\e}}
       \langle
      \A \nabla |v_{\e}|, \de_{h}\nabla |v_{\e}|\rangle\, dx
       }
       $$
       where $E_{\e}=\{x\in \O\ |\ |u(x)|>\e\}$. Then
       $$
       \displaylines{
       \int_{\O} \langle \A\nabla u, \nabla [\r^{2}_{\e}(|u|)\,
       u]\rangle\, dx =
       \int_{\O}\langle \A\nabla v_{\e}, \nabla v_{\e}\rangle\, dx -\cr
       (1-2/p) \int_{\O}{1\over |v_{\e}|} \left(
        \langle v_{\e}\A \nabla |v_{\e}|, \nabla v_{\e}\rangle-
        \langle \A\nabla v_{\e},v_{\e}\, \nabla |v_{\e}|\rangle\right)\, dx
-\cr
       ( 1-2/p)^{2} \int_{\O}
      \langle \A
       \nabla |v_{\e}|, \nabla |v_{\e}|\rangle\, dx +R(\e)
       }
       $$
       where
       $$\display{R(\e)=
       (1-2/p) \int_{\O\setminus E_{\e}}{1\over |v_{\e}|} \left(
        v_{\e}\langle \A  \nabla v_{\e}|,\nabla v_{\e}\rangle-
        \langle \A \nabla v_{\e},v_{\e}\, \nabla |v_{\e}|\rangle\right)\,
dx -\cr
       ( 1-2/p)^{2} \int_{\O\setminus E_{\e}}
       \langle \A
       \nabla |v_{\e}|, \nabla |v_{\e}|\rangle\, dx.
       }
       $$

       It is proved in \cite{langer} that if $u\in
       C^{2}(\overline{\O})$ and $u|_{\de\O}=0$, then
       \begin{equation}
	   \lim_{\e\to 0}\e^{r}\int_{\O\setminus E_{\e}}|\nabla u|^{2}dx =0
	   \label{eq:lang1}
       \end{equation}
       for any $r>-1$.
       Since
       $$
       |\, \nabla|v_{\e}|\,| = \left|\Re
       \left({\overline{v}_{\e}\nabla v_{\e}\over |v_{\e}|}\,
       \chi_{E_{0}}\right)\right|\leq |\nabla v_{\e}|= \e^{p-2\over
       2}|\nabla u|
       $$
       in $E_{0}\setminus E_{\e}$, we obtain
       $$
       \left| \int_{\O\setminus E_{\e}} \langle \A\nabla|v_{\e}|,
       \nabla|v_{\e}|\rangle\, dx \right| \leq K\, \e^{p-2}
       \int_{\O\setminus E_{\e}}
       |\nabla u|^{2} dx \to 0 
       $$
       as $\e\to 0$.
       We have also
       $$
      |v_{\e}|^{-1}
       \left| \langle v_{\e}\A \nabla |v_{\e}|, \nabla v_{\e}\rangle-
        \langle \A\nabla v_{\e},v_{\e}\, \nabla |v_{\e}|\rangle\right|\leq
       K\, \e^{p-2}|\nabla u|^{2}
       $$
       and thus $R(\e)=o(1)$ as $\e\to 0$.

      We have proved that
      \beginrighe{poti0}
       \Re\int_{\O} \langle \A \nabla u, \nabla[\r^{2}_{\e}(|u|)\,
       u]\rangle\, dx =\Re\Big[
       \int_{\O} \langle\A \nabla v_{\e}, \nabla v_{\e}\rangle\,
       dx -\crg
       (1-2/p) \int_{\O}\langle(
       \A -\A^{*}) \, \nabla |v_{\e}|, |v_{\e}|^{-1}\overline{v}_{\e} \nabla v_{\e}\rangle
       dx -\crg
        (1-2/p)^{2} \int_{\O}
       \langle \A
       \nabla |v_{\e}|, \nabla |v_{\e}|\rangle\, dx\Big] + o(1).
  \endrighe

       By means of similar computations, we find by the identity
       $$\display{
       \int_{\O}\lan {\bf b}\nabla u, |u|^{p-2}u\ran dx =
       \int_{\O\setminus E_{\e}}\lan {\bf b}\nabla u, |u|^{p-2}u\ran
       dx-\cr
       (1-2/p)\int_{E_{\e}}\lan {\bf b},|v_{\e}|\nabla(|v_{\e}|)\ran
       dx + \int_{E_{\e}}\lan {\bf b}\nabla v_{\e},v_{\e}\ran dx
       }
       $$
      that
      \beginrighe{poti1}
       \Re\int_{\O}\lan {\bf b}\nabla u, |u|^{p-2}u\ran dx =
       \crg
       \int_{\O}\lan \Re ({\bf b}/p),\nabla(|v_{\e}|^{2})\ran dx -
       \int_{\O}\lan \Im{\bf b},\Im (\overline{v}_{\e}\nabla v)\ran dx
       + o(1).
      \endrighe

Moreover
\beginrighe{poti2}
\int_{\O}|u|^{p}dx = \int_{E_{\e}}|u|^{p}dx + \int_{\O\setminus
E_{\e}}|u|^{p}dx =\crg
\int_{E_{\e}}|v_{\e}|^{2}dx +\int_{\O\setminus E_{\e}}|u|^{p}dx =
\int_{\O}|v_{\e}|^{2}dx + o(1).
\endrighe

Equalities \reff{poti0}, \reff{poti1} and \reff{poti2} lead to
\beginrighe{casino}
\Re\int_{\O} (\langle \A \nabla u, \nabla[\r^{2}_{\e}(|u|)\,
       u]\rangle -
      \lan {\bf b}\nabla u, |u|^{p-2}u\ran
      -a|u|^{p})dx =\crg
      \Re\Big[
	     \int_{\O} \langle\A \nabla v_{\e}, \nabla v_{\e}\rangle\,
	     dx -\crg
	     -(1-2/p) \int_{\O}\langle(
	     \A -\A^{*}) \, \nabla |v_{\e}|, \nabla v_{\e}\rangle
	     )v_{\e} |v_{\e}|^{-1} dx -\crg
	     - (1-2/p)^{2} \int_{\O}
	     \langle \A
	     \nabla |v_{\e}|, \nabla |v_{\e}|\rangle\, dx\Big]+\crg
	     \int_{\O} \Re (\nabla^{t}({\bf b}/p) |v_{\e}|^{2} dx +
		    \int_{\O}\lan \Im{\bf b},\Im (\overline{v}_{\e}\nabla
v)\ran dx
		    -\crg
		    \int_{\O}\Re a\, |v_{\e}|^{2}dx + o(1).
\endrighe

As far as the left-hand side of \reff{casino} is concerned, we have
$$\display{
       \int_{\O} \langle \A \nabla u, \nabla [\r^{2}_{\e}(|u|)\,
       u]\rangle\, dx =\cr
       \e^{p-2}\int_{\O\setminus E_{\e}}
	\langle \A \nabla u, \nabla
       u\rangle\, dx
       + \int_{E_{\e}}
	\langle \A \nabla u, \nabla (|u|^{p-2}u)\rangle\, dx.
       }
       $$
    and then
	$$\display{
 \lim_{\e\to 0}\Re\int_{\O} (\langle \A \nabla u, \nabla[\r^{2}_{\e}(|u|)\,
       u]\rangle -
      \lan {\bf b}\nabla u, |u|^{p-2}u\ran
      -a|u|^{p})dx =\cr
	\int_{\O} \langle \nabla u, \nabla (|u|^{p-2}u)\rangle-
      \lan {\bf b}\nabla u, |u|^{p-2}u\ran
      -a|u|^{p})dx.
      }
	       $$
	       
       Letting $\e\to 0$ in 
       \reff{casino}, we
       complete the proof of the necessity.

\textit{Sufficiency.}
       Suppose that \reff{upos} holds.
       Let $v\in \Xi$ and let  $u_{\e}$ be defined by \reff{eq:defge}.
       We have $u_{\e}\in \Xi$ and arguing as in the necessity
       part of Lemma \ref{lemma:1}, we find \reff{zac1}, \reff{zac2} 
       and
       \reff{zac4}.
       These limit relations lead to \reff{vpos}  for any $v\in \Xi$ and
       thus
       \reff{vpos} is true for any $v\in \Hspt^{1}(\O)$
       (see the proof of  Lemma \ref{lemma:2}). 
       In view of Lemma \ref{lemma:2},
       the form $\elle$ is $L^{p}$-dissipative.
       \end{proof}

\begin{theorem}\label{th:ancoraequiv}
    The operator $A$ is $L^{p}$-dissipative if and only if the form
    $\elle$ is $L^{p}$-dissipative.
\end{theorem}

\begin{proof}

    \textit{Necessity.}
    Let $u\in \Xi$ and $g_{\e}=(|u|^{2}+\e^{2})^{1\over 2}$.
    Since $g_{\e}^{p-2}\overline{u} \in \Xi$ we have
    $$
    -\int_{\O}\langle \nabla^t(\A\nabla u), u\rangle g_{\e}^{p-2}dx =
    \int_{\O}\langle \A \nabla u, \nabla(g_{\e}^{p-2}u)\rangle dx
    $$
    and since
    $$
    \de_{h}(g_{\e}^{p-2}\overline{u}) = (p-2) g_{\e}^{p-4}\Re (\langle
    \de_{h}u,u\rangle)\, \overline{u} + g_{\e}^{p-2}\de_{h}\overline{u}
    $$
    we have also
    $$\display{
    \de_{h}(g_{\e}^{p-2}\overline{u}) = \cr
    \cases{
    (p-2)|u|^{p-4}\Re (\langle
    \de_{h}u,u\rangle)\, \overline{u} + |u|^{p-2}\de_{h}\overline{u}
    = \de_{h}(|u|^{p-2}\overline{u}) & if $x\in F_{0}$\cr
    \e^{p-2}\de_{h}\overline{u} & if $x\in \O\setminus F_{0}$.\cr}
    }
    $$

  We find,  keeping in mind \reff{eq:lang1}, that
    $$
    \lim_{\e\to 0}\int_{\O}\langle  \A \nabla u,
    \nabla(g_{\e}^{p-2}u)\rangle dx=
    \int_{\O}
    \langle \A\nabla u,\nabla (|u|^{p-2}u) \rangle dx\, .
    $$

    On the other hand, using  Lemma 3.3 in \cite{langermazya}, we see
    that
    $$
    \lim_{\e\to 0}\int_{\O}\langle \nabla^t(\A\nabla u), u\rangle
g_{\e}^{p-2}dx=
    \int_{\O}\langle \nabla^t(\A\nabla u), u\rangle |u|^{p-2}dx.
    $$

    Then
    \beginrighe{ultima}
    \grande
    -\int_{\O}\langle \nabla^t(\A\nabla u), u\rangle |u|^{p-2}dx =
    \grande \int_{\O}
    \langle \A \nabla u,\nabla(|u|^{p-2}u) \rangle dx
    \endrighe
    for any $u\in \Xi$. Hence
    $$\display{
    -\int_{\O}\lan Au,u\ran |u|^{p-2}dx = \cr
    \int_{\O}(\lan \A \nabla u,\nabla(|u|^{p-2}u)\ran -
	   \lan {\bf b}\nabla u, |u|^{p-2}u\ran -a\, |u|^{p})dx \, .
	   }
   $$

   Therefore \reff{upos} holds. We can conclude now that
   the form $\elle$ is
   $L^{p}$-dissipative, because of Lemma \ref{lemma:equiv}.

    \textit{Sufficiency.} Given $u\in \Dom(A)$, we
    can find a sequence $\{u_{n}\}\subset \Xi$ such that $u_{n}\to u$ in
    $W^{2,p}(\O)$. Keeping in mind \reff{ultima}, we have
    $$\display{
    -\int_{\O}\langle Au,u\rangle |u|^{p-2}dx = -
    \lim_{n\to\infty}\int_{\O}\langle Au_{n},u_{n}\rangle |u_{n}|^{p-2}dx =\cr
    \lim_{n\to\infty}\int_{\O}\lan \A
    \nabla u_{n},\nabla(|u_{n}|^{p-2}u_{n})\ran
    - \lan {\bf b}\nabla u_{n}, |u_{n}|^{p-2}u_{n}\ran -a\, |u_{n}|^{p})dx.
    }
    $$

   Since $\elle$ is $L^{p}$-dissipative, \reff{upos} holds for any
   $u\in \Xi$ and \reff{eq:defclass} is true for any $u\in\Dom(A)$.
\end{proof}

\begin{definition}   We say that the operator $A$ is
$L^{p}$-quasi-dissipative if there exists $\o\geq 0$ such that  $A-\o
I$ is $L^{p}$-dissipative, i.e.
$$
    \Re \int_{\O}\lan Au,u\ran |u|^{p-2}dx \leq \o\, \Vert
    u\Vert_{p}^{p}
$$
for any $u\in\Dom(A)$.
\end{definition}

\begin{lemma}
    The operator \reff{eq:Acompleto} is
    $L^{p}$-quasi-dissipative if and only if there exists $\o\geq 0$ such
that
    \beginrighe{eq:22new}
	\grande
	    \Re \int_{\O}\Big[ \lan \A\nabla v,\nabla v\ran -
	    (1-2/p)\lan (\A-\A^{*})\nabla(|v|),|v|^{-1}\overline{v}
	    \nabla v\ran -\\
	    \grande
	    (1-2/p)^{2}\lan \A \nabla(|v|),\nabla(|v|)\ran
    \Big] dx
    +
 \int_{\O}\lan\Im {\bf b},\Im(\overline{v}\nabla
    v)\ran\, dx  +
    \cr
    \grande
    \int_{\O}\Re(\dive ({\bf b}/p)- a
    )|v|^{2} dx
	    \geq -\o \int_{\O}|v|^{2}dx
	\endrighe
	for any $v\in\Hspt^{1}(\O)$.
\end{lemma}

\begin{proof}
   The result follows  from Lemma \ref{lemma:2}.
\end{proof}

The next result permits to determine the best interval of $p$'s for 
   which the operator 
   \begin{equation}
       Au=\nabla^{t}(\A\nabla u)
       \label{eq:Acorto}
   \end{equation}
   is $L^{p}$-dissipative. 
   We set
   $$
   \l=\inf_{(\xi,x)\in {\cal M}}{\lan \Re\A(x)\xi,\xi\ran \over
   |\lan \Im\A(x)\xi,\xi\ran|}
   $$
   where ${\cal M}$ is the set of $(\xi,x)$ with $\xi\in\R^{n}$, 
   $x\in\O$ such that
   $\lan \Im\A(x)\xi,\xi\ran\neq 0$.

   \begin{corollary}\label{cor:5}
       Let $A$ be the 
       operator \reff{eq:Acorto}.
       Let us suppose that the matrix $\Im\A$ is symmetric and that
       \begin{equation}
	   \lan \Re\A(x)\xi,\xi\ran \geq 0
	   \label{eq:newposit}
       \end{equation}
       for any $x\in\O$, $\xi\in\R^{n}$. 
       If $\Im\A(x)=0$  for any $x\in\O$, $A$ is $L^{p}$-dissipative
       for any $p>1$.
   If $\Im\A$ does not vanish identically on $\O$, 
       $A$ is $L^{p}$-dissipative if and only if
       \begin{equation}\label{eq:bestp}
	   2+2\l(\l-\sqrt{\l^{2}+1}) \leq p \leq 
	   2+2\l(\l+\sqrt{\l^{2}+1}) .
   \end{equation}
    \end{corollary}
    
   \begin{proof}

When $\Im\A(x)=0$ for any $x\in\O$, the statement follows  from 
       Theorem \ref{th:1new}. Let us assume that $\Im\A$ does not 
       vanish identically;
      note that this implies ${\cal M}\neq \emptyset$.
      
       \textit{Necessity.} If the operator \reff{eq:Acorto}
       is $L^{p}$-dissipative, 
       Theorem \ref{th:1new} shows that
       \begin{equation}
	  |p-2|\, |\lan \Im\A(x)\xi,\xi\ran| \leq 2 \sqrt{p-1}\,
	  \lan \Re\A(x)\xi,\xi\ran
	  \label{eq:w0ter}
      \end{equation}
       for any $x\in\O$, $\xi\in\R^{n}$.
       In particular we have
       $$
       {|p-2|\over 2\sqrt{p-1}} \leq {\lan \Re\A(x)\xi,\xi\ran \over
   |\lan \Im\A(x)\xi,\xi\ran|}
       $$
       for any $(\xi,x)\in{\cal M}$ and then
       $$
       {|p-2|\over 2\sqrt{p-1}} \leq \l.
       $$
       
This inequality is equivalent to \reff{eq:bestp}.  

      \textit{Sufficiency.} If \reff{eq:bestp} holds, we have
      $(p-2)^{2}\leq 4(p-1)\l^{2}$. Note that $p>1$, because
      $2+2\l(\l-\sqrt{\l^{2}+1}) >1$.
      
      Since $\l\geq 0$ in view of \reff{eq:newposit}, we find
      $|p-2|\leq 2\sqrt{p-1}\l$ and \reff{eq:w0ter} is true
      for any $(\xi,x)\in{\cal M}$. On the other hand,
      if $x\in\O$ and $\xi\in\R^{n}$ with
      $(\xi,x)\notin{\cal M}$, \reff{eq:w0ter} is trivially satisfied and
      then it holds for any $x\in\O$, $\xi\in\R^{n}$. 
      Theorem \ref{th:1new} gives the result.
   \end{proof}

   The next Corollary  provides a characterization of operators which are
   $L^{p}$-dissipative only for $p=2$.
   
   \begin{corollary}
       Let $A$ be as in Corollary \ref{cor:5}. The operator $A$ is
       $L^{p}$-dissipative only for $p=2$ if and only if
	   $\Im\A$ does not  vanish identically and $\l=0$.
   \end{corollary}

   \begin{proof}
       Inequalities \reff{eq:bestp} are satisfied only for $p=2$ if and only if
       $\l(\l-\sqrt{\l^{2}-1})=\l(\l+\sqrt{\l^{2}+1})$ and this happens if 
       and only if $\l=0$.
       Thus the result is a consequence of  Corollary \ref{cor:5}.
   \end{proof}
   
   From now on we suppose that the operator 
      is strongly elliptic in 
      $\O$ in
      the sense that
  $$
	  \lan\Re\A(x)\xi,\xi\ran > 0
  $$
      for any $x\in\overline{\O}$, $\xi\in\R^{n}\setminus \{0\}$.

      We have proved that, if $\Im\A$ is symmetric,
      the algebraic condition \reff{eq:24} is necessary 
      and sufficient for the $L^{p}$-dissipativity of 
      the operator \reff{eq:Acorto}. We have shown that this is
      not true for the more general operator \reff{eq:Acompleto}.
      The next result shows that condition \reff{eq:24} is 
      necessary and sufficient for the $L^{p}$-quasi-dissipativity
      of \reff{eq:Acompleto}.
 We emphasize that here we do not require the symmetry of $\Im\A$.

\begin{theorem}\label{th:quasi}
    The strongly elliptic operator \reff{eq:Acompleto}
    is $L^{p}$-quasi-dissipative if and only if
    \begin{equation}
	|p-2|\, |\lan \Im\A(x)\xi,\xi\ran| \leq 2 \sqrt{p-1}\, \lan
\Re\A(x)\xi,\xi\ran
        \label{eq:w0}
    \end{equation}
   for any $x\in\O$, $\xi\in\R^{n}$.
\end{theorem}

\begin{proof}

   \textit{Necessity.} By using the functions $X$, $Y$ introduced in
   Corollary \ref{cor:-1}, we write condition \reff{eq:22new} in the
   form
$$\display{
   \int_{\O}\Big\{ {4\over p\,p'}
   \lan \Re \A X,X\ran + \lan \Re\A Y,Y\ran
	    +\cr
	     2 \lan (p^{-1}\Im\A + p'^{-1}\Im\A^{*}) X,Y\ran
	    +\lan \Im{\bf b},Y\ran |v| +\cr
	    \Re\left[\dive( {\bf b}/p) - a+\o
	    \right]|v|^{2}\Big\} dx  \geq 0\, .
 }
 $$

   As in the proof of Corollary \ref{cor:-1}, this inequality
  implies 
  $$\display{
	      {4\over p\, p'}\,
		  \int_{\O}\lan \Re\A \nabla\r,\nabla\r\ran dx  +
		  \m^{2}\int_{\O}\lan \Re\A\nabla\r,\nabla\r\ran  dx  +
		  \crg
		  2\m \int_{\O}
	   \lan (p^{-1}\Im\A + p'^{-1}\Im\A^{*})
		   \nabla\r, \nabla\r\ran dx
		   +\cr
	   \m\int_{\O}\r\lan \Im{\bf b}, \nabla\r \ran dx +
		   \int_{\O} \Re\left[\dive\left( {\bf b}/p \right) - a +\o
		      \right]\r^{2} dx
			      \geq 0
		  }
  $$
	for any $\r\in \Cspt^{1}(\O)$, $\m\in\R$.
	Since
	$$
	\lan\Im\A^{*}\nabla\r,\nabla\r\ran = -\lan
	\Im\A^{t}\nabla\r,\nabla\r\ran = -\lan \Im\A\nabla\r,\nabla\r\ran
	$$
   we have
   $$\display{
		 {4\over p\, p'}\,
		     \int_{\O}\lan \Re\A \nabla\r,\nabla\r\ran dx  +
		     \m^{2}\int_{\O}\lan \Re\A\nabla\r,\nabla\r\ran  dx -
		     \crg
		     2(1-2/p)\m \int_{\O}
	      \lan \Im\A
		      \nabla\r, \nabla\r\ran dx
		      +\cr
	      \m\int_{\O}\r\lan \Im{\bf b}, \nabla\r \ran dx +
		      \int_{\O} \Re\left[\dive\left( {\bf b}/p \right) - a +\o
			 \right]\r^{2} dx
				 \geq 0
		    }
		    $$
		      for any $\r\in \Cspt^{1}(\O)$, $\m\in\R$.

Taking $\r(x)=\p(x)\cos\lan\xi,x\ran$ and $\r(x)=\p(x)\sin\lan\xi,x\ran$ 
with
$\p\in\Cspt^{1}(\O)$ and arguing as in the proof of
Corollary \ref{cor:-1}, we find
$$\display{
\int_{\O}\lan \B\nabla\p, \nabla\p\ran dx + \int_{\O}\lan
\B\xi,\xi\ran \p^{2}dx
+ \cr
\m\int_{\O}\lan \Im {\bf b},\nabla\p\ran \p\, dx +
\int_{\O}\Re\left[\dive\left( {\bf b}/p \right) - a +\o
			 \right]\p^{2}dx \geq 0\, ,
}
			 $$
where $\m\in\R$ and
$$
\B={4\over p\,p'}\Re\A + \m^{2}\Re\A -2(1-2/p)\m \Im\A\, .
$$

Because of the arbitrariness of $\xi$ we see that
$$
\int_{\O}\lan
\B\xi,\xi\ran \p^{2}dx \geq 0
$$
for any $\p\in\Cspt^{1}(\O)$. Hence
$\lan \B\xi,\xi\ran \geq 0$, i.e.
$$
{4\over p\,p'}\lan\Re\A\xi,\xi\ran  + \m^{2}\lan\Re\A\xi,\xi\ran
-2(1-2/p)\m \lan\Im\A\xi,\xi\ran \geq 0
$$
for any $x\in\O$, $\xi\in \R^{n}$, $\m\in\R$. Inequality \reff{eq:w0}
follows from the
arbitrariness of $\m$.

   \textit{Sufficiency.} Assume first that $\Im\A$ is
   symmetric.
  By repeating  the first part of the proof of sufficiency of Theorem
  \ref{th:1new},
 we find that  \reff{eq:w0} implies
  \begin{equation}
      {4\over p\,p'}\lan \Re\A\xi,\xi\ran + \lan\Re\A\h,\h\ran
      -2(1-p/2)\lan\Im\A\xi,\h\ran \geq 0
      \label{eq:w2}
  \end{equation}
  for any $x\in\O$, $\xi,\h\in\R^{n}$.

  In order to prove \reff{eq:22new}, it is not
  restrictive to suppose
  $$
  \Re(\dive ({\bf b}/p)- a)=0.
   $$

   Since $A$ is strongly elliptic, there exists a non singular real matrix $\Cm \in
   C^{1}(\overline{\O})$ such that
   $$
   \lan\Re\A\h,\h\ran = \lan \Cm \h, \Cm \h\ran
   $$
   for any $\h\in\R^{n}$. Setting
   $$
   \Sm = (1-2/p) (\Cm^{t})^{-1}\Im\A,
   $$
   we have
   $$
   |\Cm \h-\Sm \xi|^{2}= \lan \Re \A\h,\h\ran -2(1-p/2)\lan
   \Im\A\xi,\h\ran +|\Sm \xi|^{2}.
   $$

   This leads to the identity
   \beginrighe{identity}
   {4\over p\,p'}\lan \Re\A\xi,\xi\ran + \lan\Re\A\h,\h\ran
	 -2(1-p/2)\lan\Im\A\xi,\h\ran =\crg
	 |\Cm \h-\Sm \xi|^{2} + {4\over p\,p'}\lan \Re\A\xi,\xi\ran -|\Sm
\xi|^{2}
	 \endrighe
 for any $\xi,\h\in\R^{n}$.
   In view of \reff{eq:w2}, putting  $\h=\Cm^{-1}\Sm \xi$ in \reff{identity},
   we obtain
   \begin{equation}
       {4\over p\,p'}\lan \Re\A\xi,\xi\ran -|\Sm \xi|^{2} \geq 0
       \label{diseg}
   \end{equation}
   for any $\xi\in\R^{n}$.

   On the other hand, we may write
   $$\display{
   \lan \Im{\bf b}, Y\ran = \lan (\Cm^{-1})^{t}\Im{\bf b}, \Cm Y\ran =\cr
   \lan (\Cm^{-1})^{t}\Im{\bf b}, \Cm Y-\Sm X\ran +  \lan
(\Cm^{-1})^{t}\Im{\bf b},
   \Sm X\ran \, .
   }
   $$

   By  the Cauchy inequality 
   $$\display{
   \int_{\O}\lan (\Cm^{-1})^{t}\Im{\bf b}, \Cm Y- \Sm X\ran |v|\, dx \geq \cr
   -\int_{\O}|\Cm Y-\Sm X|^{2}dx - {1\over 4}\int_{\O}|(\Cm^{-1})^{t}\Im{\bf
   b}|^{2}|v|^{2}dx
   }
   $$
   and, integrating by parts,
   $$\display{
   \int_{\O}\lan (\Cm^{-1})^{t}\Im{\bf b},
   \Sm X\ran |v|\, dx = {1\over 2}\int_{\O}\lan (\Cm^{-1}\Sm )^{t}\Im{\bf b},
   \nabla(|v|^{2})\ran\, dx =\cr
   -{1\over 2}\int_{\O} \nabla^{t}((\Cm^{-1}\Sm )^{t}\Im{\bf b})
      \, |v|^{2}\, dx\, .
      }
   $$

   This implies that there exists $\o \geq 0$ such that
$$
       \int_{\O}\lan \Im{\bf b}, Y\ran \,|v|\, dx \geq
       -\int_{\O}|\Cm Y-\Sm X|^{2}dx - \o\int_{\O}|v|^{2}dx
$$
   and then, in view of \reff{identity},
   $$\display{
   \int_{\O}\Big\{ {4\over p\,p'}
     \lan \Re \A X,X\ran + \lan \Re\A Y,Y\ran
	      +\cr
	      \grande 2 (1-p/2)\lan \Im\A  X,Y\ran
	      +\lan \Im{\bf b},Y\ran |v| \Big\}dx \geq\cr
	\int_{\O}\left({4\over p\,p'}\lan \Re\A X,X\ran -|\Sm X|^{2}\right)dx
	-\o\int_{\O}|v|^{2}dx\ .
}
	      $$

	      Inequality \reff{diseg} gives the result.

	      We have proved the sufficiency under the assumption
$\Im\A^{t}=\Im\A$.
	     In the general case, the operator $A$ can be written in the form
	     $$
	     Au= \nabla^{t}((\A+\A^{t})\nabla u)/2 + {\bf c}\nabla u + a u
	     $$
	     where
	     $$
	     {\bf c}= \nabla^{t}(\A-\A^{t})/2 +{\bf b}.
	     $$

	     Since $(\A+\A^{t})$ is symmetric, we know that $A$ is
	     $L^{p}$-quasi-dissipative if and only if
	     $$
	     |p-2|\, |\lan \Im(\A+\A^{t})\xi,\xi\ran| \leq 2 \sqrt{p-1}\, \lan
	     \Re(\A+\A^{t})\xi,\xi\ran
	     $$
	     for any $\xi\in\R^{n}$, which is exactly condition \reff{eq:w0}.
\end{proof}

\begin{corollary}
    Let $A$ be the strongly elliptic operator \reff{eq:Acompleto}. 
    If $\Im\A(x)=0$  for any $x\in\O$, $A$ is $L^{p}$-quasi-dissipative
	   for any $p>1$.
       If $\Im\A$ does not vanish identically on $\O$, 
	   $A$ is $L^{p}$-quasi-dissipative if and only if
	   \reff{eq:bestp} holds.
\end{corollary}

\begin{proof}
    The proof is similar to that 
 of Corollary \ref{cor:5}, the role of Theorem  \ref{th:1new}
being played by
    Theorem \ref{th:quasi}. 
\end{proof}

    We give 
    a  criterion for the $L^{p}$-contractivity of the semigroup
    generated by $A$.

    \begin{theorem}\label{th:contr}
    	Let $A$ be the  strongly elliptic operator \reff{eq:Acorto} 
    with
    $\Im\A=\Im\A^{t}$. 
    The operator $A$ generates a contraction semigroup
   on $L^{p}$ if and only if
     \begin{equation}
	|p-2|\, |\lan \Im\A(x)\xi,\xi\ran| \leq 2 \sqrt{p-1}\,
	\lan \Re\A(x)\xi,\xi\ran
        \label{eq:w0bis}
    \end{equation}
   for any $x\in\O$, $\xi\in\R^{n}$. 
    \end{theorem}

    \begin{proof}
	
	\textit{Sufficiency.}
	It is a classical result that the operator $A$ defined on 
	\reff{eq:DomA}
	and acting in $L^{p}(\O)$ is a densely defined closed operator
	(see \cite{agmon}, 
	    \cite[Theorem 1, p.302]{mazyashaposh}).
  
    From Theorem \ref{th:1new} we know that the form $\elle$ is
    $L^{p}$-dissipative and Theorem \ref{th:ancoraequiv} shows that
    $A$ is $L^{p}$-dissipative.
    Finally the formal adjoint operator
    $$
    A^{*}u=\nabla^{t}(\A^{*}\nabla u)
    $$
    with $\Dom(A^{*})=W^{2,p'}(\O)\cap \Wspt^{1,p'}(\O)$, is the
    adjoint operator of $A$  and since
    $\Im\A^{*}=\Im(\A^{*})^{t}$ and  \reff{eq:w0bis} can be written as
    \begin{equation}
        |p'-2|\, |\lan \Im\A^{*}(x)\xi,\xi\ran| \leq 2 \sqrt{p'-1}\,
            \lan \Re\A^{*}(x)\xi,\xi\ran ,
        \label{eq:condp'}
    \end{equation}
	we have
    also the $L^{p'}$-dissipativity of $A^{*}$.

    The result is a consequence of the following well known result: 
    \textit{if $A$ is a
densely defined
    closed operator and
    if both $A$ and $A^{*}$ are dissipative, then $A$ is the
    infinitesimal generator of a $C_{0}$ contraction semigroup}
     (see, e.g.,
    \cite[p.15]{pazy}).
    
    \textit{Necessity.} If $A$ generates a contraction semigroup on 
    $L^{p}$, 
    it is $L^{p}$-dissipa\-tive. Therefore \reff{eq:w0bis} holds because of 
    Theorem \ref{th:1new}.
    \end{proof}

    Let us assume that either $A$ has lower order terms
    or they are absent and $\Im\A$ is not symmetric. The next Theorem 
    gives a criterion for the $L^{p}$-quasi-contractivity 
     of the semigroup generated by $A$
    (i.e. the $L^{p}$-contractivity of the semigroup generated by $A-\o I$).
 
    \begin{theorem}
	Let $A$ be the strongly elliptic operator \reff{eq:Acompleto}. 
       The operator $A$ generates a 
       quasi-contraction semigroup
   on $L^{p}$ if and only if
    \reff{eq:w0bis} holds
   for any $x\in\O$, $\xi\in\R^{n}$. 
    \end{theorem}
    
    \begin{proof}
	
	\textit{Sufficiency.} Let us consider $A$ as an operator defined on 
	\reff{eq:DomA} and acting in $L^{p}(\O)$. 
	As in the proof of Theorem \ref{th:contr},
	one can see that $A$ is a densely defined closed operator and 
	that the 
	formal adjoint coincides with the adjoint $A^{*}$. Theorem 
	\ref{th:quasi} shows that $A$ 
	is $L^{p}$-quasi-dissipative. On the other hand,  condition \reff{eq:condp'}
	holds and then $A^{*}$ is $L^{p'}$-quasi-dissipative.
	As in Theorem \ref{th:contr}, this implies that $A$ generates a 
	quasi-contraction semigroup on $L^{p}$.
	
	\textit{Necessity.} If $A$ generates a quasi-contraction semigroup 
	on $L^{p}$, $A$ is $L^{p}$-quasi-dissipative and \reff{eq:w0bis} 
	holds.
	\end{proof}

\end{document}